\documentclass[reqno]{amsproc}
\usepackage[latin1]{inputenc}
\usepackage{amssymb}
\usepackage{hyperref}
\usepackage{amsmath}
\usepackage{latexsym}
\usepackage{cite}
\usepackage{amssymb}
\usepackage{amsfonts}
\usepackage{amscd}
\usepackage{xcolor}
\usepackage{amstext,amsmath,amssymb,amsfonts}
\newtheorem{theorem}{Theorem}

\newtheorem{definition}[theorem]{Definition}
\newtheorem{example}[theorem]{Example}

\newtheorem{proposition}[theorem]{Proposition}
\newtheorem{remark}[theorem]{Remark}

\newcommand{\A}{\mathcal{A}}

\newcommand{\B}{\mathcal{B}}

\newcommand{\M}{\mathcal{M}}

\newcommand{\beq}{\begin{eqnarray}}
\newcommand{\eeq}{\end{eqnarray}}
\newcommand{\beqs}{\begin{eqnarray*}}
\newcommand{\eeqs}{\end{eqnarray*}}
\newcommand{\bpro}{\begin{pro}}
\newcommand{\epro}{\end{pro}}
\newcommand{\blem}{\begin{lem}}
\newcommand{\elem}{\end{lem}}
\newcommand{\bdfn}{\begin{dfn}}
\newcommand{\edfn}{\end{dfn}}
\newcommand{\bcor}{\begin{cor}}
\newcommand{\ecor}{\end{cor}}
\newcommand{\bthm}{\begin{thm}}
\newcommand{\ethm}{\end{thm}}
\newcommand{\bex}{\begin{ex}}
\newcommand{\eex}{\end{ex}}
\newcommand{\brmk}{\begin{rmk}}
\newcommand{\ermk}{\end{rmk}}
\newcommand{\bpr}{\begin{pr}}
\newcommand{\epr}{\end{pr}}
\newcommand{\benum}{\begin{enumerate}}
\newcommand{\eenum}{\end{enumerate}}
\newcommand{\bitem}{\begin{itemize}}
\newcommand{\eitem}{\end{itemize}}

\chardef\bslash=`\\
\numberwithin{equation}{section}
\numberwithin{table}{section}
\numberwithin{theorem}{section}
\DeclareMathOperator{\id}{id}

\title[ Hom-coassociative ternary coalgebras and infinitesimal bialgebras]{ Hom-coassociative ternary  coalgebras and infinitesimal bialgebras  \footnote{Preprint: ICMPA-MPA/2015/08 } }

\author{Mahouton Norbert Hounkonnou$^\ast$}
\address[$\ast$]{University of Abomey-Calavi,
International Chair in Mathematical Physics and Applications,
ICMPA-UNESCO Chair, 072 BP 50, Cotonou, Rep. of Benin}
\email{norbert.hounkonnou@cipma.uac.bj, with copy to hounkonnou@yahoo.fr}
\author{Gb\^ev\`ewou Damien  Houndedji$^\dagger$}
\address[$\dagger$]{University of Abomey-Calavi,
International Chair in Mathematical Physics and Applications,
ICMPA-UNESCO Chair, 072 BP 50, Cotonou, Rep. of Benin}
\email{ houndedjid@gmail.com}
\begin{document}
\maketitle

\today

\bigskip
\begin{abstract}
 In this work, the partially and totally hom-coassociative ternary coalgebras are constructed and discussed.   Their {infinitesimal} bialgebraic structures
 are also investigated. The related dual space structures and their properties are elucidated.  
\end{abstract}

\section{Introduction}
The hom-algebraic  structures  were originated from  quasi-deformations of Lie algebras of vector fields. 
The latters  gave rise to   quasi-Lie algebras,  defined as  generalized Lie structures in which the 
skew-symmetry and Jacobi conditions are twisted. The first examples were concerned with $q$-deformations 
of the Witt and Virasoro algebras (see \cite{[chaichian1]}, \cite{[chaichian2]}, \cite{[chaichian3]}, 
\cite{[chaichian4]}, \cite{[curtright]}, \cite{[daskaloyannis]}, \cite{[hu]}, \cite{[kassel]}, 
\cite{[liu]}). Motivated by {those works} and  examples { applying} the 
general quasi-deformation construction provided by   Larsson and co-workers (\cite{[hartwig]}, 
\cite{[larsson1]}, \cite{[larsson2]}), and the possibility of   studying,  within the same framework, 
such well-known generalizations of Lie algebras as the color and super Lie algebras, the quasi-Lie 
algebras, subclasses of quasi-hom-Lie algebras and hom-Lie algebras were introduced in 
(\cite{[hartwig]}, \cite{[larsson1]}, \cite{[larsson2]}, \cite{[larsson3]}). It is worth noticing here 
that, in the subclass of hom-Lie algebras, the  skew-symmetry is untwisted, whereas the Jacobi identity 
is twisted by a single linear map and contains three terms as in Lie algebras, reducing to ordinary Lie 
algebras when the twisting linear map is the identity map.

The  hom-associative algebras {are} generalizations of  the  associative algebras, where the 
associativity law is twisted by a linear map. In (\cite{[makhlouf]}),  it was shown in particular that 
the commutator product, defined using the multiplication in a hom-associative algebra, naturally leads  
to a hom-Lie algebra. The hom-Lie-admissible algebras and more general $G$-hom-
associative algebras with subclasses of hom-Vinberg and pre-hom-Lie algebras, generalizing twisted  
Lie-admissible algebras, $G$-associative algebras, Vinberg and pre-Lie algebras {were also  introduced}.  For these 
classes of algebras, the operation of taking commutator leads to hom-Lie algebras as well. The enveloping 
algebras of hom-Lie algebras were discussed in (\cite{[yau1]}).
The fundamentals of the formal deformation theory and associated cohomology structures for hom-Lie 
algebras were recently considered in (\cite{[makhlouf2]}). Simultaneously,  Yau  developed elements of 
homology for hom-Lie algebras in (\cite{[yau2]}). In (\cite{[makhlouf1]}) and (\cite{[makhlouf3]}), 
Ataguema {\it et al } elaborated the theory of hom-coalgebras and related structures. They introduced 
hom-coalgebra structure, leading to the notions of hom-bialgebras and hom-Hopf algebras, proved some 
fundamental properties and provided examples. They also defined the concepts of hom-Lie admissible hom-
coalgebra generalizing the admissible coalgebra introduced in \cite{[goze]}, and {performed} their classification based on subgroups of the {underlying} symmetric group.

The present  work aims at  defining and discussing  properties of   partially and totally hom-
coassociative ternary coalgebras.  The {infinitesimal} bialgebraic structures of  hom-associative ternary 
algebras and  their relation with corresponding hom-associative ternary algebras are investigated. The 
dual structures of such categories of algebras  are also studied.  

For the clarity of our presentation, let us conclude this section by briefly recalling main definitions 
and some relevant known results, which can be useful in the sequel. In this context, we  mainly focuss on 
the definitions of totally  and partially hom-associative ternary algebras given in \cite{[P1]}, 
generalizing the classical totally  and partially associative ternary algebras.  Then, we provide  a 
method of their construction (see also  \cite{[P1]}) starting from a totally  or partially associative algebra 
and their algebra endomorphism.
\subsection{Basic definitions}
 We give the definitions of totally hom-associative algebras and partially hom-associative algebras.
  Let $\mathcal{K}$ be a commutative field of characteristics zero. 
\begin{definition}
A totally hom-associative ternary algebra is  a $\mathcal{K}$-vector space $\mathcal{T}$, {endowed with}  a trilinear operation $\mu: \mathcal{T}\otimes\mathcal{T}\otimes\mathcal{T}\rightarrow \mathcal{T}$,  and a familly $\lbrace \alpha_{1},  \alpha_{2}\rbrace$ of linear maps $\alpha_{1}, \alpha_{2}: \mathcal{T}\rightarrow \mathcal{T}$ satisfying, for all $x_{1}, x_{2}, x_{3}, x_{4}, x_{5}\in \mathcal{A},$
\beq \label{qt1}
\mu(\mu(x_{1}\otimes x_{2}\otimes x_{3})\otimes \alpha_{1}(x_{4})\otimes \alpha_{2}(x_{5}))=\cr \mu(\alpha_{1}(x_{1})\otimes \mu(x_{2}\otimes x_{3}\otimes x_{4})\otimes \alpha_{2}(x_{5}))=\cr \mu(\alpha_{1}(x_{1})\otimes \alpha_{2}(x_{2})\otimes \mu(x_{3}\otimes x_{4}\otimes x_{5})).
\eeq 
\end{definition}

\begin{definition}
A partially hom-associative ternary algebra is  a $\mathcal{K}$-vector space $\mathcal{P}$, {equipped with}  a trilinear operation $\mu: \mathcal{P}\otimes\mathcal{P}\otimes\mathcal{P}\rightarrow \mathcal{P}$,  and a familly $\lbrace \alpha_{1},  \alpha_{2}\rbrace$ of linear maps $\alpha_{1}, \alpha_{2}: \mathcal{P}\rightarrow \mathcal{P}$ satisfying for all $x_{1}, x_{2}, x_{3}, x_{4}, x_{5}\in \mathcal{P}$
\beq \label{qp1}
&&\mu(\mu(x_{1}\otimes x_{2}\otimes x_{3})\otimes \alpha_{1}(x_{4})\otimes \alpha_{2}(x_{5}))+ \mu(\alpha_{1}(x_{1})\otimes \mu(x_{2}\otimes x_{3}\otimes x_{4})\otimes \alpha_{2}(x_{5}))+\cr
 &&\mu(\alpha_{1}(x_{1})\otimes \alpha_{2}(x_{2})\otimes \mu(x_{3}\otimes x_{4}\otimes x_{5}))=0.
\eeq 
\end{definition}

\begin{remark} 
A weak totally hom-associative ternary algebra is given by the identity 
\beqs
\mu(\mu(x_{1}\otimes x_{2}\otimes x_{3})\otimes \alpha_{1}(x_{4})\otimes \alpha_{2}(x_{5}))= \mu(\alpha_{1}(x_{1})\otimes \alpha_{2}(x_{2})\otimes \mu(x_{3}\otimes x_{4}\otimes x_{5})).
\eeqs  
\end{remark}
{
\begin{definition}
A hom-associative ternary algebra $(\A, \mu, \alpha_{1}, \alpha_{2})$ is said multiplicative if
\beq
\alpha_{i}\circ\mu= \mu\circ\alpha_{i}^{\otimes^{3}}, \ i=1, 2.
\eeq
\end{definition}
}
The morphisms of hom-associative ternary algebras are defined as follows.
\begin{definition}
Let $(\mathcal{A}, \mu, \alpha)$ and $(\mathcal{A}', \mu', \alpha')$ be two totally hom-associative   (resp. partially hom-associative) ternary algebras, where $\alpha=(\alpha_{i})_{i=1, 2}$ and $\alpha'=(\alpha'_{i})_{i=1, 2}$. A linear map $\rho: \mathcal{A}\rightarrow \mathcal{A}'$ is a totally hom-associative ternary algebras (resp. partially hom-associative ternary algebras) morphism if it satisfies
\beqs
\rho(\mu(x_{1}, x_{2}, x_{3}))= \mu'(\rho(x_{1}), \rho(x_{2}), \rho(x_{3})) \mbox{ and } \rho\circ \alpha_{i}=\alpha'_{i}\circ \rho \ \ \forall i= 1, 2.
\eeqs
\end{definition}
\subsection{Construction of hom-associative ternary algebras}
The following theorem \cite{[P1]} provides a way to construct a totally (resp. a partially)  hom-associative algebra starting from a totally  (resp. partially) associative ternary algebra and ternary algebra endomorphism.

\begin{theorem}\label{teo1}
Let $(\mathcal{A}, \mu)$ be a totally  (resp. partially)  associative ternary algebra 
and let $\rho: \mathcal{A}\rightarrow \mathcal{A}$ be a totally  (resp. partially) hom-
associative ternary algebra endomorphism. We set $\mu_{\rho}= \rho\circ\mu$ and $\tilde{\rho}= 
(\rho_{i})_{i=1,..., 5}$ with $\rho_{i}= \rho$. Then $(\mathcal{A}, \mu_{\rho}, \tilde{\rho})$ is a 
totally (resp. partially) hom-associative ternary algebra.
Moreover, suppose that $(\mathcal{A}', \mu')$ is another totally  (resp. 
partially) associative ternary algebra and $\rho': \mathcal{A}'\rightarrow \mathcal{A}'$ be a totally 
 (resp. partially)  hom-associative ternary algebra endomorphism. If $f: \mathcal{A}\rightarrow 
\mathcal{A}'$ is a totally  (resp. partially)  associative ternary 
algebra morphism that is $ f\circ \rho= \rho'\circ f,$ then 
\beqs
f: (\mathcal{A}, \mu_{\rho}, \tilde{\rho})\rightarrow (\mathcal{A}', \mu'_{\rho'}, \tilde{\rho'})
\eeqs
is a totally  (resp. partially) hom-associative ternary algebra morphism.
\end{theorem}
Therefore the totally associativity ternary identity implies the totally hom-associativity ternary identity for the  ternary algebras.
\begin{example}\label{ep1}
Let $\mathcal{P}$ be a $2$-dimensional vector space  with a basis $\lbrace e_{1}, e_{2} \rbrace$. The trilinear product $\mu$ on $\mathcal{P}$ defined by
\beqs
\mu(e_{1}\otimes e_{1}\otimes e_{1})= e_{2}; \ \  \mu(e_{i}\otimes e_{j}\otimes e_{k})=0, i, j, k=1, 2 \mbox{ with } (i, j, k)\neq (1, 1, 1)
\eeqs
defines a partially associatve ternary algebra.
Its automorphisms are given with respect to the same basis by the matrices of the form: $ \left(
\begin{array}{cc}
a & 0  \\
b & a^{3}  \\
\end{array}
\right),$ with $a, b\in \mathcal{K}$ and $a\neq 0.$\\
We deduce, using Theorem \ref{teo1}, the $2-$dimensional partially hom-associative ternary algebra, $(\mathcal{P}, \tilde{\mu}, \rho),$ given, with respect to the basis $\lbrace e_{1}, e_{2} \rbrace$, by the following nontrivial product:
\beqs
\tilde{\mu}(e_{1}, e_{1}, e_{1})= a^{3}e_{2},
\eeqs
and the linear map:
\beqs
\rho(e_{1})&=& a e_{1} + b e_{2}\cr
\rho(e_{2})&=& a^{3} e_{2}.
\eeqs 
\end{example}
\begin{example}\label{et1}
Let $\mathcal{T}$ be a $2$-dimensionnal vector space with a basis $\lbrace e_{1}, e_{2}\rbrace$. Consider the following product, defined by
\beqs
&& \mu: \mathcal{T}\otimes \mathcal{T}\otimes\mathcal{T}\rightarrow \mathcal{T},\cr
&&\mu(e_{1}\otimes e_{1}\otimes e_{1})= e_{1} \ \ \ \mu(e_{2}\otimes e_{2}\otimes e_{1})= e_{1} + e_{2}\cr
&&\mu(e_{1}\otimes e_{1}\otimes e_{2})= e_{2} \ \ \ \mu(e_{2}\otimes e_{2}\otimes e_{2})= e_{1} + 2e_{2}\cr
&&\mu(e_{1}\otimes e_{2}\otimes e_{1})= e_{2} \ \ \ \mu(e_{1}\otimes e_{2}\otimes e_{2})= e_{1} + e_{2}\cr
&&\mu(e_{2}\otimes e_{1}\otimes e_{1})= e_{2} \ \ \ \mu(e_{2}\otimes e_{1}\otimes e_{2})= e_{1} + e_{2},
\eeqs
defines a totally associatve ternary algebra.
The following matrices are automorphisms of this totally associative ternary algebra:
\beqs
\left(
\begin{array}{cc}
1 & 0  \\
0 & 1  \\
\end{array}
\right), \ \
\left( 
\begin{array}{cc}
-1 & 0  \\
0 & -1  \\
\end{array}
\right), \ \
\left( 
\begin{array}{cc}
-1 & -1  \\
0 & 1  \\
\end{array}
\right), \cr
\left(
\begin{array}{cc}
-\dfrac{1}{\sqrt{5}} & -\dfrac{3}{\sqrt{5}}  \\
\dfrac{2}{\sqrt{5}} & \dfrac{1}{\sqrt{5}}  \\
\end{array}
\right), \ \
\left( 
\begin{array}{cc}
1 & 1  \\
0 & -1  \\
\end{array}
\right), \ \
\left( 
\begin{array}{cc}
\dfrac{1}{\sqrt{5}} & \dfrac{3}{\sqrt{5}}  \\
-\dfrac{2}{\sqrt{5}} & -\dfrac{1}{\sqrt{5}}  \\
\end{array}
\right).    
\eeqs
We deduce, for example. using Theorem \ref{teo1}, the following $2$-dimensional totally hom-associative ternary algebras $(\mathcal{T}, \tilde{\mu}_{i}, \rho_{i})$ for $i= 1, 2,$ which are not totally associative. They are given, with respect to a basis $\lbrace e_{1}, e_{2} \rbrace$, by the following product:
\beqs
&&\tilde{\mu}_{1}(e_{1}\otimes e_{1}\otimes e_{1})= e_{1} \ \ \ \tilde{\mu}_{1}(e_{2}\otimes e_{2}\otimes e_{1})= 2e_{1} - e_{2}\cr
&&\tilde{\mu}_{1}(e_{1}\otimes e_{1}\otimes e_{2})= e_{1}-e_{2} \ \ \ \tilde{\mu}_{1}(e_{2}\otimes e_{2}\otimes e_{2})= 3e_{1} - 2e_{2}\cr
&&\tilde{\mu}_{1}(e_{1}\otimes e_{2}\otimes e_{2})= 2e_{1}-e_{2} \ \ \ \tilde{\mu}_{1}(e_{1}\otimes e_{2}\otimes e_{1})= e_{1} - e_{2}\cr
&&\tilde{\mu}_{1}(e_{2}\otimes e_{1}\otimes e_{1})= e_{1}-e_{2} \ \ \ \tilde{\mu}_{1}(e_{2}\otimes e_{1}\otimes e_{2})= 2e_{1} - e_{2}
\eeqs
and the linear map
\beqs
\rho_{1}(e_{1})&=& e_{1}\cr
\rho_{1}(e_{2})&=& e_{1}-e_{2}.
\eeqs 
A second structure is given by 
\beqs
&&\tilde{\mu}_{2}(e_{1}\otimes e_{1}\otimes e_{1})= \dfrac{1}{\sqrt{5}}e_{1}- \dfrac{2}{\sqrt{5}}e_{2} \ \ \ \tilde{\mu}_{2}(e_{2}\otimes e_{2}\otimes e_{1})= \dfrac{4}{\sqrt{5}}e_{1} - \dfrac{3}{\sqrt{5}}e_{2}\cr
&&\tilde{\mu}_{2}(e_{1}\otimes e_{1}\otimes e_{2})= \dfrac{3}{\sqrt{5}}e_{1}-\dfrac{1}{\sqrt{5}}e_{2} \ \ \ \tilde{\mu}_{2}(e_{2}\otimes e_{2}\otimes e_{2})= \dfrac{7}{\sqrt{5}}e_{1} - \dfrac{4}{\sqrt{5}}e_{2}\cr
&&\tilde{\mu}_{2}(e_{1}\otimes e_{2}\otimes e_{2})= \dfrac{4}{\sqrt{5}}e_{1}-\dfrac{3}{\sqrt{5}}e_{2} \ \ \ \tilde{\mu}_{2}(e_{1}\otimes e_{2}\otimes e_{1})= \dfrac{3}{\sqrt{5}}e_{1} - \dfrac{1}{\sqrt{5}}e_{2}\cr
&&\tilde{\mu}_{2}(e_{2}\otimes e_{1}\otimes e_{1})= \dfrac{3}{\sqrt{5}}e_{1}-\dfrac{1}{\sqrt{5}}e_{2} \ \ \ \tilde{\mu}_{2}(e_{2}\otimes e_{1}\otimes e_{2})= \dfrac{4}{\sqrt{5}}e_{1} - \dfrac{3}{\sqrt{5}}e_{2}
\eeqs
and the linear map
\beqs
\rho_{2}(e_{1})&=& \dfrac{1}{\sqrt{5}}e_{1}-\dfrac{2}{\sqrt{5}}e_{2}\cr
\rho_{2}(e_{2})&=& \dfrac{3}{\sqrt{5}}e_{1}-\dfrac{1}{\sqrt{5}}e_{2}.
\eeqs 
\end{example}
Finally,  for a linear map $\alpha: \mathcal{A}\rightarrow \mathcal{A},$ we denote the dual (linear) map by $\alpha^{\ast}:\mathcal{A}^{\ast}\rightarrow \mathcal{A}^{\ast}$ given by
\beq \label{eqalpha}
\langle \alpha^{\ast}(f), x\rangle= \langle f, \alpha(x)\rangle \ f\in \mathcal{A}^{\ast},\ x\in \mathcal{A}. 
\eeq
\section{Hom-coassociative ternary coalgebras}
In this section, we introduce and develop the concepts of totally and partially hom-coassociative  ternary coalgebras. 

Let $(\mathcal{A}, \mu, (\alpha_{i})_{i=1, 2})$ be a hom-associative  ternary algebra, and $\mathcal{A}^{\ast},$ its dual space. Then, we get the dual mapping $\mu^{\ast}: \mathcal{A}^{\ast}\rightarrow \mathcal{A}^{\ast}\otimes\mathcal{A}^{\ast}\otimes\mathcal{A}^{\ast}$ of $\mu$, for every $x, y, z\in \mathcal{A}$ and $\xi, \eta, \gamma \in \mathcal{A}^{\ast}$, 
\beq \label{qdual}
\langle \mu^{\ast}(\xi), x\otimes y\otimes z\rangle= \langle\xi, \mu(x, y, z) \rangle,
\eeq 
\beq
\langle \xi\otimes\eta\otimes\gamma, x\otimes y\otimes z\rangle= \langle\xi, x\rangle \langle\eta, y\rangle \langle\gamma, z\rangle,
\eeq
where $\langle, \rangle$ is the natural nondegenerate symmetric bilinear form on the vector space $\mathcal{A}\oplus \mathcal{A}^{\ast},$ defined by $\langle \xi, x\rangle= \xi(x), \xi\in \mathcal{A}^{\ast}, x\in \mathcal{A}$.
By the definition of hom-associative ternary algebras, we have $Im(\mu^{\ast})\subseteq \mathcal{A}^{\ast}\otimes\mathcal{A}^{\ast}\otimes\mathcal{A}^{\ast},$ and
\beqs
(\mu^{\ast}\otimes \alpha^{\ast}_{1}\otimes \alpha^{\ast}_{2})\circ\mu^{\ast}= (\alpha^{\ast}_{1}\otimes \mu^{\ast}\otimes \alpha^{\ast}_{2})\circ\mu^{\ast}= (\alpha^{\ast}_{1}\otimes \alpha^{\ast}_{2}\otimes \mu^{\ast})\circ\mu^{\ast} \mbox { and } 
\eeqs 
\beqs
(\mu^{\ast}\otimes \alpha^{\ast}_{1}\otimes \alpha^{\ast}_{2} +\alpha^{\ast}_{1}\otimes \mu^{\ast}\otimes \alpha^{\ast}_{2}+ \alpha^{\ast}_{1}\otimes \alpha^{\ast}_{2}\otimes \mu^{\ast})\circ\mu^{\ast}=0, 
\eeqs 
for the totally hom-associative ternary algebras and the partially hom-associative ternary algebras, respectively. That is, for every $x_{1}, x_{2}, x_{3}, x_{4}, x_{5}\in \mathcal{A}$ and $\xi\in \mathcal{A}^{\ast}:$
\beqs
&&\langle(\mu^{\ast}\otimes\alpha^{\ast}_{1}\otimes\alpha^{\ast}_{2})\circ\mu^{\ast}(\xi), x_{1}\otimes x_{2}\otimes x_{3}\otimes x_{4}\otimes x_{5}\rangle\cr
&&=\langle\xi, \mu\circ(\mu\otimes\alpha_{1}\otimes\alpha_{2})(x_{1}\otimes x_{2}\otimes x_{3}\otimes x_{4}\otimes x_{5})\rangle\cr
&&=\langle\xi, \mu(\mu(x_{1}\otimes x_{2}\otimes x_{3})\otimes \alpha_{1}(x_{4})\otimes \alpha_{2}(x_{5}))\rangle \cr
&&=\langle\xi, \mu(\alpha_{1}(x_{1})\otimes \mu(x_{2}\otimes x_{3}\otimes x_{4})\otimes \alpha_{2}(x_{5}))\rangle\cr
&&=\langle(\alpha^{\ast}_{1}\otimes\mu^{\ast}\otimes\alpha^{\ast}_{2})\circ\mu^{\ast}(\xi), x_{1}\otimes x_{2}\otimes x_{3}\otimes x_{4}\otimes x_{5}\rangle\cr
&&=\langle\xi, \mu(\alpha_{1}(x_{1})\otimes \alpha_{2}(x_{2})\otimes \mu(x_{3}\otimes x_{4}\otimes x_{5}))\rangle\cr
&&=\langle(\alpha^{\ast}_{1}\otimes\alpha^{\ast}_{2}\otimes\mu^{\ast})\circ\mu^{\ast}(\xi), x_{1}\otimes x_{2}\otimes x_{3}\otimes x_{4}\otimes x_{5}\rangle\cr 
\eeqs 
and
\beqs
&&\langle(\mu^{\ast}\otimes \alpha^{\ast}_{1}\otimes\alpha^{\ast}_{2} + \alpha^{\ast}_{1}\otimes \mu^{\ast}\otimes \alpha^{\ast}_{2} + \alpha^{\ast}_{1}\otimes \alpha^{\ast}_{2}\otimes \mu^{\ast})\circ\mu^{\ast}(\xi), x_{1}\otimes x_{2}\otimes x_{3}\otimes x_{4}\otimes x_{5}\rangle\cr
=&&\langle\xi, \mu\circ(\mu\otimes \alpha_{1}\otimes \alpha_{2} +\alpha_{1}\otimes \mu\otimes \alpha_{2} +\alpha_{1}\otimes \alpha_{2}\otimes \mu)(x_{1}\otimes x_{2}\otimes x_{3}\otimes x_{4}\otimes x_{5})\rangle,
\eeqs
respectively.
 \begin{definition}\label{THCTC}
 A totally hom-coassociative ternary coalgebra $(\mathcal{T}, \Delta, (\alpha_{i})_{i=1, 2})$ is a $\mathcal{K}$-vector space $\mathcal{T}$ with a linear mapping $\Delta: \mathcal{T}\rightarrow  \mathcal{T}\otimes\mathcal{T}\otimes\mathcal{T}$ and a familly $\lbrace \alpha_{1},  \alpha_{2}\rbrace$ of linear maps $\alpha_{1}, \alpha_{2}: \mathcal{T}\rightarrow \mathcal{T}$ satisfying 
 \beq 
 (\Delta\otimes \alpha_{1}\otimes \alpha_{2})\circ\Delta= (\alpha_{1}\otimes \Delta\otimes \alpha_{2})\circ\Delta= (\alpha_{1}\otimes \alpha_{2}\otimes \Delta)\circ\Delta. 
 \eeq
  \end{definition}
  
 \begin{definition}\label{WTHCTC}
 A weak totally hom-coassociative ternary coalgebra $(\mathcal{W}, \Delta, (\alpha_{i})_{i=1, 2})$ is a $\mathcal{K}$-vector space $\mathcal{W}$ with a linear mapping $\Delta: \mathcal{W}\rightarrow  \mathcal{W}\otimes\mathcal{W}\otimes\mathcal{W}$ and a familly $\lbrace \alpha_{1},  \alpha_{2}\rbrace$ of linear maps $\alpha_{1}, \alpha_{2}: \mathcal{W}\rightarrow \mathcal{W}$ satisfying
 \beq 
 (\Delta\otimes \alpha_{1}\otimes \alpha_{2})\circ\Delta= (\alpha_{1}\otimes \alpha_{2}\otimes \Delta)\circ\Delta.
 \eeq
  \end{definition}  
  
 \begin{definition}\label{PHCTC}
  A partially hom-coassociative ternary coalgebra $(\mathcal{P}, \Delta, (\alpha_{i})_{i=1, 2})$ is a $\mathcal{K}$-vector space $\mathcal{P}$ with a linear mapping $\Delta: \mathcal{P}\rightarrow  \mathcal{P}\otimes\mathcal{P}\otimes\mathcal{P}$ and a familly $\lbrace \alpha_{1},  \alpha_{2}\rbrace$ of linear maps $\alpha_{1}, \alpha_{2}: \mathcal{P}\rightarrow \mathcal{P}$ satisfying
 \beq 
 (\Delta\otimes \alpha_{1}\otimes \alpha_{2} + \alpha_{1}\otimes \Delta\otimes \alpha_{2} + \alpha_{1}\otimes \alpha_{2}\otimes \Delta)\circ\Delta= 0. 
 \eeq
\end{definition}
{
\begin{definition}
A hom-coassociative ternary coalgebra $(C, \Delta, (\alpha_{i})_{i=1, 2})$ is said comultiplicative if
\beq
\alpha^{\otimes^{3}}_{i}\circ\Delta= \Delta\circ\alpha_{i}, \ i=1, 2.
\eeq
\end{definition}
}
Given the above definitions, we  now study partially hom-coassociative ternary coalgebras in terms   of structural constants. Let $(\mathcal{P}, \Delta, (\alpha_{i})_{i=1, 2})$ be a partially hom-coassociative ternary coalgebra with a basis $e_{1}, ....., e_{n}$. Assume that 
\beq\label{q1}
&&\Delta(e_{l})= \sum_{1\leq r, s, t\leq n}c^{l}_{rst}e_{r}\otimes e_{s}\otimes e_{t}, \ \ c^{l}_{rst}\in \mathcal{K}, 1\leq l\leq n \mbox{ and } \cr
&&\alpha_{i}(e_{j})=\sum^{n}_{k=1}y^{ij}_{k}e_{k}, \ j=1,...,n.
\eeq
Then, we have
\beqs
&&(\Delta\otimes\alpha_{1}\otimes\alpha_{2} + \alpha_{1}\otimes\Delta\otimes\alpha_{2} + \alpha_{1}\otimes\alpha_{2}\otimes\Delta)\circ\Delta(e_{l})\cr
&&=(\Delta\otimes\alpha_{1}\otimes\alpha_{2} + \alpha_{1}\otimes\Delta\otimes\alpha_{2} + \alpha_{1}\otimes\alpha_{2}\otimes\Delta)\left( \sum_{1\leq r, s, t\leq n}c^{l}_{rst}e_{r}\otimes e_{s}\otimes e_{t}\right) \cr
&&=\sum_{1\leq r, s, t\leq n}c^{l}_{rst}\left( \Delta(e_{r})\otimes \alpha_{1}(e_{s})\otimes \alpha_{2}(e_{t}) + \alpha_{1}(e_{r})\otimes \Delta(e_{s})\otimes \alpha_{2}(e_{t}) + \alpha_{1}(e_{r})\otimes \alpha_{2}(e_{s})\otimes \Delta(e_{t})\right) \cr
&&=\sum_{1\leq r, s, t\leq n}c^{l}_{rst} \lbrace \left(\sum_{1\leq i, j, k\leq n}c^{r}_{ijk}e_{i}\otimes e_{j}\otimes e_{k}\right) \otimes\left(\sum^{n}_{q=1}y^{1s}_{q}e_{q}\right) \otimes \left( \sum^{n}_{p=1}y^{2t}_{p}e_{p}\right)  +\cr
&& \left(\sum^{n}_{q=1}y^{1r}_{q}e_{q}\right) \otimes \left(\sum_{1\leq i, j, k\leq n}c^{s}_{ijk} e_{i}\otimes e_{j}\otimes e_{k}\right) \otimes \left(\sum^{n}_{p=1}y^{2t}_{p}e_{p}\right)  + \cr 
&&\left(\sum^{n}_{q=1}y^{1r}_{q}e_{q}\right)\otimes
 \left(\sum^{n}_{p=1}y^{2s}_{p}e_{p}\right)\otimes \left(\sum_{1\leq i, j, k\leq n}c^{t}_{ijk}e_{i}\otimes e_{j}\otimes e_{k}\right)  \rbrace  \cr
&&=\sum_{1\leq r, s, t, i, j, k, q, p\leq n}\lbrace c^{l}_{rst}c^{r}_{ijk}y^{1s}_{q}y^{2t}_{p}e_{i}\otimes e_{j}\otimes e_{k}\otimes e_{q}\otimes e_{p} +c^{l}_{rst}c^{s}_{ijk}y^{1r}_{q}y^{2t}_{p}e_{q}\otimes e_{i}\otimes e_{j}\otimes e_{k}\otimes e_{p}+ \cr && c^{l}_{rst}c^{t}_{ijk}y^{1r}_{q}y^{2s}_{p}e_{q}\otimes e_{p}\otimes e_{i}\otimes e_{j}\otimes e_{k}\rbrace \cr
&&=\sum_{1\leq r, s, t, i, j, k, q, p\leq n} \left[ c^{l}_{rst}c^{r}_{ijk}y^{1s}_{q}y^{2t}_{p} +c^{l}_{rst}c^{s}_{jkq}y^{1r}_{i}y^{2t}_{p}+  c^{l}_{rst}c^{t}_{kqp}y^{1r}_{i}y^{2s}_{j}\right] e_{i}\otimes e_{j}\otimes e_{k}\otimes e_{q}\otimes e_{p}
\eeqs 
Therefore, 
\beq \label{q2}
\sum^{n}_{r=1}\left[c^{l}_{rst}c^{r}_{ijk}y^{1s}_{q}y^{2t}_{p} +c^{l}_{rst}c^{s}_{jkq}y^{1r}_{i}y^{2t}_{p}+  c^{l}_{rst}c^{t}_{kqp}y^{1r}_{i}y^{2s}_{j}\right]=0, \ \ 1\leq i, j, k, s, t, l, p,q\leq n.
\eeq
By a similar discussion, for a totally hom-coassociative ternary coalgebra, we obtain 
\beq \label{q3}
\sum^{n}_{r=1}c^{l}_{rst}c^{r}_{ijk}y^{1s}_{q}y^{2t}_{p}= \sum^{n}_{r=1}c^{l}_{rst}c^{s}_{jkq}y^{1r}_{i}y^{2t}_{p}=\sum^{n}_{r=1} c^{l}_{rst}c^{t}_{kqp}y^{1r}_{i}y^{2s}_{j}, \ \ 1\leq i, j, k, s, t, l, p, q\leq n,
\eeq
and for a weak totally hom-coassociative ternary coalgebra:
\beq \label{q4}
\sum^{n}_{r=1}c^{l}_{rst}c^{r}_{ijk}y^{1s}_{q}y^{2t}_{p}= \sum^{n}_{r=1} c^{l}_{rst}c^{t}_{kqp}y^{1r}_{i}y^{2s}_{j}, \ \ 1\leq i, j, k, s, t, l, p, q\leq n.
\eeq
Therefore, we can set the following:
 \begin{theorem}\label{Constant structures theorem for HATC}
 Let $\mathcal{A}$ be an $n$-dimensional vector space with a basis $e_{1},....., e_{n},$
 
  $\Delta: \mathcal{A}\rightarrow \mathcal{A}\otimes \mathcal{A}\otimes \mathcal{A},\ \alpha_{i=1, 2}: \mathcal{A}\rightarrow \mathcal{A}$ be defined as (\ref{q1}). Then,
 \begin{enumerate}
 \item $(\mathcal{A}, \Delta)$ is a partially hom-coassociative ternary coalgebra if and only if the constants $c^{l}_{ijk}$and $y^{qs}_{p},$  $1\leq i, j, k, p, s\leq n;\ q=1, 2$ satisfy the identity (\ref{q2});
 \item $(\mathcal{A}, \Delta)$ is a totally hom-coassociative ternary coalgebra if and only if the constants $c^{l}_{ijk}$and $y^{qs}_{p},$  $1\leq i, j, k, p, s\leq n;\ q=1, 2$ satisfy  the identity (\ref{q3});
 \item $(\mathcal{A}, \Delta)$ is a weak totally hom-coassociative ternary coalgebra if and only if the constants $c^{l}_{ijk}$and $y^{qs}_{p},$  $1\leq i, j, k, p, s\leq n;\ q=1, 2$ satisfy the identity  (\ref{q4}).
 \end{enumerate}
  \end{theorem}
  
  Now, let $(\mathcal{A}, \mu, (\alpha_{i})_{i=1, 2})$ be a partially hom-associative ternary algebra with a basis $e_{1}, e_{2},....., e_{n}.$ The multiplication $\mu$ of $\mathcal{A}$ in this basis is {given} as follows:
  
\beq
\mu(e_{r}, e_{s}, e_{t})=\sum^{n}_{l=1}c^{l}_{rst}e_{l}, \ \  c^{l}_{rst}\in \mathcal{K}, \ \ 1\leq r, s, t\leq n.
\eeq  
Using the condition (\ref{qp1}), we have:
\beqs
&&\mu(\mu(e_{r}, e_{s}, e_{t}), \alpha_{1}(e_{i}), \alpha_{2}(e_{j})) + \mu(\alpha_{1}(e_{r}), \mu(e_{s}, e_{t}, e_{i}), \alpha_{2}(e_{j})) + \mu(\alpha_{1}(e_{r}), \alpha_{2}(e_{s}), \mu(e_{t}, e_{i}, e_{j}))\cr
&&=\mu\left( \sum^{n}_{l=1}c^{l}_{rst} e_{l}, \sum^{n}_{q=1}y^{1i}_{q} e_{q}, \sum^{n}_{p=1}y^{2j}_{p} e_{p}\right)  + \mu\left( \sum^{n}_{q=1}y^{1r}_{q} e_{q}, \sum^{n}_{l=1}c^{l}_{sti} e_{l}, \sum^{n}_{p=1}y^{2j}_{p} e_{p}\right)   +\cr &&\mu\left( \sum^{n}_{q=1}y^{1r}_{q} e_{q}, \sum^{n}_{p=1}y^{2s}_{p} e_{q}, \sum^{n}_{l=1}c^{l}_{tij} e_{l}\right) \cr
&&=\sum^{n}_{l, p, q=1}\left( c^{l}_{rst}y^{1i}_{q}y^{2j}_{p} \mu(e_{l}, e_{q}, e_{p})+ c^{l}_{sti}y^{1r}_{q}y^{2j}_{p} \mu(e_{q}, e_{l}, e_{p}) + c^{l}_{tij}y^{1r}_{q}y^{2s}_{p} \mu(e_{q}, e_{p}, e_{l})\right)  \cr
&&= \sum^{n}_{l, p, q, k=1}\left( c^{l}_{rst}c^{k}_{lpq}y^{1i}_{q}y^{2j}_{p} + c^{l}_{sti}c^{k}_{qlp}y^{1r}_{q}y^{2j}_{p} + c^{l}_{tij}c^{k}_{qpl}y^{1r}_{q}y^{2s}_{p}\right)e_{k}=0.
\eeqs
  Therefore, $\sum^{n}_{l=1}\left(c^{l}_{rst}c^{k}_{lpq}y^{1i}_{q}y^{2j}_{p} + c^{l}_{sti}c^{k}_{qlp}y^{1r}_{q}y^{2j}_{p} + c^{l}_{tij}c^{k}_{qpl}y^{1r}_{q}y^{2s}_{p}\right)=0,$ that is, \\$ \lbrace c^{l}_{i_{1}i_{2}i_{3}, 1\leq i_{1}, i_{2},
    i_{3}\leq n};\ y^{ij}_{k}, 1\leq k, j\leq n, i=1, 2 \rbrace$ satisfies the identity (\ref{q2}).
    
   A similar discussion for a totally hom-associative ternary algebra and a weak totally hom-associative ternary algebra yields, respectively:
        
         $\sum^{n}_{l=1}c^{l}_{rst}c^{k}_{lpq}y^{1i}_{q}y^{2j}_{p}= \sum^{n}_{l=1}c^{l}_{sti}c^{k}_{qlp}y^{1r}_{q}y^{2j}_{p}= 
   \sum^{n}_{l=1}c^{l}_{tij}c^{k}_{qpl}y^{1r}_{q}y^{2s}_{p},$ that is, $ \lbrace c^{l}_{i_{1}i_{2}i_{3}, 1\leq i_{1}, i_{2}, i_{3}\leq n};\\ y^{ij}_{k}, 1\leq k, j\leq n, i=1, 2 \rbrace$ satisfies the identity (\ref{q3}),
   
   and
   
   $\sum^{n}_{l=1}c^{l}_{rst}c^{k}_{lpq}y^{1i}_{q}y^{2j}_{p}= \sum^{n}_{l=1}c^{l}_{tij}c^{k}_{qpl}y^{1r}_{q}y^{2s}_{p},$ that is, $ \lbrace c^{l}_{i_{1}, i_{2}, i_{3}, 1\leq i_{1}, i_{2}, i_{3}\leq n};\ y^{ij}_{k}, 1\leq k, j\leq n,\\ i=1, 2 \rbrace$ satisfies the identity (\ref{q4}).
   
   Let $\mathcal{A}^{\ast}$ be the dual space of partially hom-associative ternary algebra
$(\mathcal{A}, \mu, (\alpha_{i})_{i=1, 2})$, and $e^{\ast}_{1},....., e^{\ast}_{n}$ be the dual
 basis of $e_{1},....., e_{n}, \langle e^{\ast}_{1}, e_{j}\rangle= \delta_{ij}, 1\leq i, j\leq.$
 Assume that\\ $\mu^{\ast}: \mathcal{A}^{\ast}$ $ \rightarrow$ $ \mathcal{A}^{\ast}\otimes
\mathcal{A}^{\ast}\otimes \mathcal{A}^{\ast}$ and $\alpha^{\ast}_{i=1, 2}: \mathcal{A}^{\ast}\rightarrow \mathcal{A}^{\ast}$  are the dual mapping of $\mu,$ and $\alpha_{i=1, 2},$ {are} defined by (\ref{qdual}) and (\ref{eqalpha}), respectively. Then, for every $1\leq l\leq n,$ we get
    \beq
    \mu^{\ast}(e^{\ast}_{l})= \sum_{1\leq r, s, t\leq n}c^{l}_{rst}e^{\ast}_{r}\otimes e^{\ast}_{s}\otimes e^{\ast}_{t}.
    \eeq
    According to the identities (\ref{q1}) and (\ref{q2}), ($\mathcal{A}^{\ast}, \mu^{\ast}, (\alpha^{\ast}_{i})_{i=1, 2}$) is a partially hom-coassociative ternary coalgebra.  

Conversely, if $(\mathcal{A}, \Delta, (\alpha_{i})_{i=1, 2})$ is a partially hom-coassociative ternary coalgebra with a basis $e_{1}, ....., e_{n}$ satisfying (\ref{q1}), $\mathcal{A}^{\ast}$ is the dual space of $\mathcal{A}$ with the dual basis $ e^{\ast}_{1}, ...., e^{\ast}_{n}.$ Then, the dual mapping $\Delta^{\ast}:$ $\mathcal{A}^{\ast}$ $ \rightarrow$ $ \mathcal{A}^{\ast}\otimes \mathcal{A}^{\ast}\otimes \mathcal{A}^{\ast}$ of $\Delta$ satisfies, for every $\xi, \eta, \gamma \in \mathcal{A}^{\ast}, x\in \mathcal{A},$
\beq \label{eqco}
\langle \Delta^{\ast}(\xi, \eta, \gamma), x\rangle= \langle\xi\otimes\eta\otimes\gamma, \Delta(x) \rangle.
\eeq 
Then, $\Delta^{\ast}(e^{\ast}_{r}, e^{\ast}_{s}, e^{\ast}_{t})= \sum^{n}_{l=1}c^{l}_{rst}e^{\ast}_{l}, c^{l}_{rst}\in \mathcal{K}, 1\leq r, s, t, l\leq n$ and  $(\Delta^{\ast}, (\alpha^{\ast}_{i})_{i=1, 2}))$ satisfies identity (\ref{q2}).
\begin{remark}
Analogous discussions remain valid for  totally  and   weak totally hom-associative ternary algebras. 
\end{remark}

Therefore, the following results are in order.
\begin{theorem} \label{teo-co-asso}
Let $\mathcal{A}$ be a vector space over a field $\mathcal{K}$, $\Delta: \mathcal{A}\rightarrow \mathcal{A}\otimes \mathcal{A}\otimes \mathcal{A}$ and a familly $\lbrace \alpha_{1},  \alpha_{2}\rbrace$ of linear maps $\alpha_{1}, \alpha_{2}: \mathcal{A}\rightarrow \mathcal{A}.$ Then,
\begin{enumerate}
\item $(\mathcal{A}, \Delta, (\alpha_{i})_{i=1, 2})$ is a partially hom-coassociative ternary coalgebra if and only if\\ $(\mathcal{A}^{\ast}, \Delta^{\ast}, (\alpha^{\ast}_{i})_{i=1, 2})$ is a partially hom-associative ternary algebra.
\item $(\mathcal{A}, \Delta, (\alpha_{i})_{i=1, 2})$ is a totally hom-coassociative ternary coalgebra if and only if\\ $(\mathcal{A}^{\ast}, \Delta^{\ast}, (\alpha^{\ast}_{i})_{i=1, 2})$ is a totally hom-associative ternary algebra.
\item $(\mathcal{A}, \Delta, (\alpha_{i})_{i=1, 2})$ is a weak totally hom-coassociative ternary coalgebra if and only if\\ $(\mathcal{A}^{\ast}, \Delta^{\ast}, (\alpha^{\ast}_{i})_{i=1, 2})$ is a weak totally hom-associative ternary algebra.
\end{enumerate}
\end{theorem}

We can also give an equivalence description of (\ref{teo-co-asso}).
\begin{theorem} \label{teo-co-asso1}
Let $\mathcal{A}$ be a vector space over a field $\mathcal{K}$, $\mu: \mathcal{A}\otimes \mathcal{A}\otimes \mathcal{A} \rightarrow \mathcal{A}$ a trilinear mapping and a familly $\lbrace \alpha_{1},  \alpha_{2}\rbrace$ of linear maps $\alpha_{1}, \alpha_{2}: \mathcal{A}\rightarrow \mathcal{A}$. Then,
\begin{enumerate}
\item $(\mathcal{A}, \mu, (\alpha_{i})_{i=1, 2})$ is a partially hom-associative ternary algebra if and only if\\ $(\mathcal{A}^{\ast}, \mu^{\ast}, (\alpha^{\ast}_{i})_{i=1, 2})$ is a partially hom-coassociative tenary coalgebra.
\item $(\mathcal{A}, \mu, (\alpha_{i})_{i=1, 2})$ is a totally hom-associative ternary algebra if and only if\\ $(\mathcal{A}^{\ast}, \mu^{\ast}, (\alpha^{\ast}_{i})_{i=1, 2})$ is a totally hom-coassociative ternary coalgebra.
\item $(\mathcal{A}, \mu, (\alpha_{i})_{i=1, 2})$ is a weak totally hom-associative ternary algebra if and only if\\ $(\mathcal{A}^{\ast}, \mu^{\ast}, (\alpha^{\ast}_{i})_{i=1, 2})$ is a weak totally hom-coassociative ternary coalgebra.
\end{enumerate}
\end{theorem}

\begin{example}
Let $(\mathcal{P}, \tilde{\mu}^{\ast}, \rho^{\ast})$ be the dual of partially hom-associative ternary
 algebra $(\mathcal{P}, \tilde{\mu}, \rho)$ in Example (\ref{ep1}). The product $\mu^{\ast}$ and the
  applications $\alpha^{\ast}_{1}=\alpha^{\ast}_{2}= \rho^{\ast}$  are given, respectively, by
\beqs
\tilde{\mu}^{\ast}(e^{\ast}_{2})=a^{3} e^{\ast}_{1}\otimes e^{\ast}_{1}\otimes e^{\ast}_{1}; \ \  \tilde{\mu}(e^{\ast}_{i}\otimes e^{\ast}_{j}\otimes e^{\ast}_{k})=0, i, j, k=1, 2 \mbox{ with } (i, j, k)\neq (1, 1, 1),
\eeqs
and the linear map
\beqs
\rho^{\ast}(e^{\ast}_{1})&=& a\ e^{\ast}_{1}\cr
\rho^{\ast}(e^{\ast}_{2})&=& b\ e^{\ast}_{1} + a^{3} e^{\ast}_{2}.
\eeqs
$(\mathcal{P}^{\ast}, \mu^{\ast}, \rho^{\ast})$ is a partially hom-coassociative ternary coalgebra.
\end{example}

\begin{example}
Let $(\mathcal{T}^{\ast}, \mu^{\ast}, \rho^{\ast})$ be the dual of totally hom-associative ternary algebra $(\mathcal{T}, \mu, \rho)$ in Example (\ref{et1}). The product $\mu^{\ast}$ on $\mathcal{T}^{\ast}$ is given by
\beqs
\tilde{\mu}^{\ast}(e^{\ast}_{1})&=& e^{\ast}_{1}\otimes e^{\ast}_{1}\otimes e^{\ast}_{1} + 2e^{\ast}_{1}\otimes e^{\ast}_{2}\otimes e^{\ast}_{2} + 2e^{\ast}_{2}\otimes e^{\ast}_{2}\otimes e^{\ast}_{1} + 3e^{\ast}_{2}\otimes e^{\ast}_{2}\otimes e^{\ast}_{2} + 2e^{\ast}_{2}\otimes e^{\ast}_{1}\otimes e^{\ast}_{2} +\cr && e^{\ast}_{1}\otimes e^{\ast}_{1}\otimes e^{\ast}_{2} + e^{\ast}_{2}\otimes e^{\ast}_{1}\otimes e^{\ast}_{1} + e^{\ast}_{1}\otimes e^{\ast}_{2}\otimes e^{\ast}_{1}\cr
\tilde{\mu}^{\ast}(e^{\ast}_{2})&=& -e^{\ast}_{1}\otimes e^{\ast}_{1}\otimes e^{\ast}_{2} - e^{\ast}_{1}\otimes e^{\ast}_{2}\otimes e^{\ast}_{2} - e^{\ast}_{2}\otimes e^{\ast}_{1}\otimes e^{\ast}_{1} - e^{\ast}_{2}\otimes e^{\ast}_{2}\otimes e^{\ast}_{1} \cr
&& -2e^{\ast}_{2}\otimes e^{\ast}_{2}\otimes e^{\ast}_{2} -e^{\ast}_{1}\otimes e^{\ast}_{2}\otimes e^{\ast}_{1} -e^{\ast}_{2}\otimes e^{\ast}_{1}\otimes e^{\ast}_{2},
\eeqs
and the linear map
\beqs
\rho^{\ast}(e^{\ast}_{1})&=& \ e^{\ast}_{1} + e^{\ast}_{2}\cr
\rho(e^{\ast}_{2})&=& -e^{\ast}_{2}.
\eeqs
$(\mathcal{T}^{\ast}, \mu^{\ast}, \rho^{\ast})$ is a totally hom-coassociative ternary coalgebra.
\end{example}
{In the sequel, and as a  matter of global analysis, the terminology of hom-coassociative ternary coalgebra will stand for partially, totally and weak totally hom-coassociative ternary coalgebras}.
\begin{definition}
Let ($\mathcal{A}_{1}, \Delta_{1}, (\alpha_{i})_{i=1, 2}$) and ($\mathcal{A}_{2}, \Delta_{2}, (\alpha'_{i})_{i=1, 2}$) be two hom-coassociative ternary coalgebras. If there is a linear isomorphism $\varphi: \mathcal{A}_{1}\rightarrow \mathcal{A}_{2}$ satisfying
\beq
(\varphi\otimes\varphi\otimes\varphi)(\Delta_{1}(e))= \Delta_{2}(\varphi(e)), \mbox{ for every } e \in \mathcal{A}_{1},\ \mbox{ and } \rho\circ \alpha_{i}=\alpha'_{i}\circ \rho \ \ \forall i= 1, 2,
\eeq 
then $(\mathcal{A}_{1}, \Delta_{1}, (\alpha_{i})_{i=1, 2})$ is isomorphic to $(\mathcal{A}_{2}, \Delta_{2}, (\alpha'_{i})_{i=1, 2}),$ and $\varphi$ is called a hom-coassociative ternary coalgebra isomorphism, defined as: 
\beq
(\varphi\otimes\varphi\otimes\varphi)\sum_{i}(a_{i}\otimes b_{i}\otimes c_{i})= \sum_{i}\varphi(a_{i})\otimes\varphi(b_{i})\otimes\varphi(c_{i}).
\eeq
\end{definition}
\begin{theorem}\label{theorem of isomorphisms and duals}
Let ($\mathcal{A}_{1}, \Delta_{1}, (\alpha_{i})_{i=1, 2}$) and ($\mathcal{A}_{2}, \Delta_{2}, (\alpha'_{i})_{i=1, 2}$) be two hom-coassociative ternary
 coalgebras. Then, $\varphi: \mathcal{A}_{1}\rightarrow \mathcal{A}_{2}$ is a hom-coassociative ternary
  coalgebra isomorphism from ($\mathcal{A}_{1}, \Delta_{1}$) to ($\mathcal{A}_{2}, \Delta_{2}$) if and
   only if the dual mapping $ \varphi^{\ast}: \mathcal{A}^{\ast}_{2}\rightarrow
    \mathcal{A}^{\ast}_{1} $ is a hom-associative ternary algebra isomorphism defined from
     ($\mathcal{A}^{\ast}_{2}, \Delta^{\ast}_{2}, (\alpha'^{\ast}_{i})_{i=1, 2}$) to ($\mathcal{A}^{\ast}_{1}, \Delta^{\ast}_{1}, (\alpha^{\ast}_{i})_{i=1, 2}$),
      such that, for every $\xi \in \mathcal{A}^{\ast}_{2}, v\in \mathcal{A}_{1}, \langle  \varphi^{\ast
      }(\xi), v\rangle= \langle\xi, \varphi(v)\rangle$.
\end{theorem}
\textbf{Proof:}
{Let} ($\mathcal{A}_{1}, \Delta_{1}, (\alpha_{i})_{i=1, 2}$) and ($\mathcal{A}_{2}, \Delta_{2}, (\alpha'_{i})_{i=1, 2}$) {be} two hom-coassociative ternary
 coalgebras. Then, ($\mathcal{A}^{\ast}_{1}, \Delta^{\ast}_{1}$) and ($\mathcal{A}^{\ast}_{2},
 \Delta^{\ast}_{2}$) are two hom-associative ternary algebras. Let $\varphi: \mathcal{A}_{1}\rightarrow 
\mathcal{A}_{2}$ be a hom-coassociative ternary coalgebra isomorphism from ($\mathcal{A}_{1},
\Delta_{1}, (\alpha_{i})_{i=1, 2}$) to ($\mathcal{A}_{2}, \Delta_{2}, (\alpha'_{i})_{i=1, 2}$). Hence, the dual mapping $\varphi^{\ast}:
 \mathcal{A}^{\ast}_{2}\rightarrow \mathcal{A}^{\ast}_{1}$ is {a} linear isomorphism, and, for every
  $\xi, \eta, \gamma \in \mathcal{A}^{\ast}_{2}$, $x \in \mathcal{A}_{1},$
  \beqs
\langle \varphi^{\ast}\Delta^{\ast}_{2}(\xi, \eta, \gamma), x\rangle &=&\langle \xi\otimes\eta\otimes\gamma, \Delta_{2}(\varphi(x))\rangle= \langle \xi\otimes\eta\otimes\gamma, (\varphi\otimes \varphi\otimes\varphi)\Delta_{1}(x)\rangle\cr
&=&\langle \varphi^{\ast}(\xi)\otimes\varphi^{\ast}(\eta)\otimes\varphi^{\ast}(\gamma), \Delta_{1}(x)\rangle= \langle \Delta^{\ast}_{1}(\varphi^{\ast}(\xi), \varphi^{\ast}(\eta), \varphi^{\ast}(\gamma)), x\rangle 
\eeqs
\beqs
\langle \rho^{\ast}\circ\alpha'^{\ast}_{i}(\xi), x\rangle = \langle \xi, \alpha'_{i}\circ\rho(x)\rangle
= \langle \xi, \rho\circ\alpha(x)\rangle
= \langle \alpha^{\ast}_{i}\circ\rho^{\ast}(\xi), x\rangle, i=1, 2.
  \eeqs
Then, $\varphi^{\ast}\Delta^{\ast}_{2}(\xi, \eta, \gamma)= \Delta^{\ast}_{1}(\varphi^{\ast}(\xi), \varphi^{\ast}(\eta), \varphi^{\ast}(\gamma))$ and $\rho^{\ast}\circ\alpha'^{\ast}_{i}= \alpha^{\ast}_{i}\circ\rho^{\ast}, i=1, 2,$ proving that  $\varphi^{\ast}$ is a hom-associative ternary algebra isomorphism.
$ \hfill \square $

\section{Trimodules and matched pairs of hom-associative ternary algebras}
In this section we introduce the concept of  trimodules over hom-associative ternary algebra and give their matched pairs. Recall that a hom-module is a pair $(V, \beta)$ in which $V$ is a vector space and $\beta: V\rightarrow V$ is a linear map. We give the definition of bihom-module as follows.
\begin{definition}
A bihom-module is a triple $(V, \beta_{1}, \beta_{2})$ in which $V$ is a vector space  and $\beta_{1}, \beta_{2}: V\rightarrow V$ are linear maps.
\end{definition}
\subsection{Trimodules and matched pairs of totally hom-associative ternary  algebras}
\begin{definition}\label{trimodule of THATA}
A trimodule structure over totally hom-associative ternary algebra $(\A, \mu, \alpha_{1}, \alpha_{2})$ on a bihom-module $(V, \beta_{1}, \beta_{2})$ is defined by the following three linear multiplication mappings: 
\beqs
&&\mathcal{L}_{\mu}: \A\otimes\A\otimes V\rightarrow V, \cr
&& \mathcal{R}_{\mu}: V\otimes\A\otimes\A\rightarrow V, \cr
&& \mathcal{M}_{\mu}: \A\otimes V\otimes\A\rightarrow V
\eeqs 
 satisfying the following compatibility conditions
\beq \label{tr1}
\mathcal{L}_{\mu}(\alpha_{1}(a), \alpha_{2}(b))(\mathcal{L}_{\mu}(c, d)(v))= \mathcal{L}_{\mu}(\mu(a, b, c), \alpha_{1}(d))\beta_{2}(v)= \mathcal{L}_{\mu}(\alpha_{1}(a), \mu(b, c, d))\beta_{2}( v),
\eeq
\beq \label{tr2}
\mathcal{R}_{\mu}(\alpha_{1}(c), \alpha_{2}(d))(\mathcal{R}_{\mu}(a, b)(v))= \mathcal{R}_{\mu}(\alpha_{2}(a), \mu(b, c, d))\beta_{1}(v)= \mathcal{R}_{\mu}(\mu(a, b, c), \alpha_{2}(d))\beta_{1}(v),
\eeq 
\beq \label{tr4}
\mathcal{M}_{\mu}(\alpha_{1}(a), \alpha_{2}(d))(\mathcal{L}{\mu}(b, c)(v))= \mathcal{L}_{\mu}(\alpha_{1}(a), \alpha_{2}(b)) (\mathcal{M}_{\mu}(c, d)(v))= \M_{\mu}(\mu(a, b, c), \alpha_{2}(d))\beta_{1}(v),
\eeq
\beq \label{tr5}
\mathcal{M}_{\mu}(\alpha_{1}(a), \alpha_{2}(d))(\mathcal{R}_{\mu}(b, c)(v))= \mathcal{R}_{\mu}(\alpha_{1}(c), \alpha_{2}(d))(\mathcal{M}_{\mu}(a, b)(v))= \mathcal{M}_{\mu}(\alpha_{1}(a), \mu(b, c, d))\beta_{2}(v),
\eeq
\beq \label{tr6}
\mathcal{R}_{\mu}(\alpha_{1}(c), \alpha_{2}(d))(\mathcal{L}_{\mu}(a, b) (v))&=& \mathcal{L}_{\mu}(\alpha_{1}(a), \alpha_{2}(b))(\mathcal{R}_{\mu}(c, d)(v))\cr
&=& \mathcal{M}_{\mu}(\alpha_{1}(a), \alpha_{2}(d))(\mathcal{M}_{\mu}(b, c)(v)),
\eeq
\beq \label{tr3}
\mathcal{M}_{\mu}(\alpha_{1}(a), \alpha_{2}(z))(\mathcal{M}_{\mu}(\alpha_{1}(b), \alpha_{2}(y))(\mathcal{M}_{\mu}(\alpha_{1}(c), \alpha_{2}(x))\beta_{1}(v))=\cr
 \mathcal{M}_{\mu}(\mu(\alpha_{1}(a), \alpha_{1}(b), \alpha_{1}(c)), \mu(\alpha_{2}(x), \alpha_{2}(y), \alpha_{2}(z)))\beta_{1}(v),
\eeq
\beq \label{tr3,2}
\mathcal{M}_{\mu}(\alpha_{1}(a), \alpha_{2}(z))(\mathcal{M}_{\mu}(\alpha_{1}(b), \alpha_{2}(y))(\mathcal{M}_{\mu}(\alpha_{1}(c), \alpha_{2}(x))\beta_{2}(v))=\cr
 \mathcal{M}_{\mu}(\mu(\alpha_{1}(a), \alpha_{1}(b), \alpha_{1}(c)), \mu(\alpha_{2}(x), \alpha_{2}(y), \alpha_{2}(z)))\beta_{2}(v),
\eeq
\beq\label{tr7}
\beta_{1}(\mathcal{L}_{\mu}(a, b)v)= \mathcal{L}_{\mu}(\alpha_{1}(a), \alpha_{2}(b))\beta_{1}(v),\ \beta_{2}(\mathcal{L}_{\mu}(a, b)v)= \mathcal{L}_{\mu}(\alpha_{1}(a), \alpha_{2}(b))\beta_{2}(v)
\eeq
\beq\label{tr8}
\beta_{1}(\mathcal{M}_{\mu}(a, b)v)= \mathcal{M}_{\mu}(\alpha_{1}(a), \alpha_{2}(b))\beta_{1}(v),\
\beta_{2}(\mathcal{M}_{\mu}(a, b)v)= \mathcal{M}_{\mu}(\alpha_{1}(a), \alpha_{2}(b))\beta_{2}(v)
\eeq
\beq\label{tr9}
\beta_{1}(\mathcal{R}_{\mu}(a, b)v)= \mathcal{R}_{\mu}(\alpha_{1}(a), \alpha_{2}(b))\beta_{1}(v), 
\beta_{2}(\mathcal{R}_{\mu}(a, b)v)= \mathcal{R}_{\mu}(\alpha_{1}(a), \alpha_{2}(b))\beta_{2}(v)
\eeq
$ \forall a, b, c, d, x, y, z\in \A, v\in V$.
\end{definition}
\begin{remark}
In the case where only the Eqs.(\ref{tr1})-(\ref{tr6}) are  satisfied, we refer to the name of quasi trimodule structure instead of simply trimodule structure.
\end{remark}

\begin{proposition}
($ \mathcal{L}_{\mu}, \mathcal{M}_{\mu}, \mathcal{R}_{\mu}, \beta_{1}, \beta_{2}, V$) is a quasi trimodule of totally hom-associative ternary algebra $(\A, \mu, \alpha_{1}, \alpha_{2})$ if and only if the direct sum {$(\mathcal{A} \oplus V, \tau, \alpha_{1} + \beta_{1}, \alpha_{2} + \beta_{2})$}  of the underlying vector spaces of $ \mathcal{A} $ and $V$  is turned into a totally hom-associative ternary algebra $\tau$ given by
\beqs
\tau[(x + a), (y + b), (z + c)]= \mu(x, y, z) + \mathcal{L}_{\mu}(x, y)(c) + \mathcal{M}_{\mu}(x, z)(b) + \mathcal{R}_{\mu}(y, z)(a),
\eeqs
for all $ x, y, z \in \mathcal{A}, a, b, c \in V$. We denote it by $ \mathcal{A} \ltimes_{\mathcal{L}_{\mu}, \mathcal{M}_{\mu}, \mathcal{R}_{\mu}, \beta_{1}, \beta_{2}} V$.
\end{proposition}
\textbf{Proof}
Let $v_{1}, v_{2}, v_{3}, v_{4}, v_{5}\in V$ and $x_{1}, x_{2}, x_{3}, x_{4}, x_{5}\in \mathcal{A}$.
Set
\beqs
&&\tau[\tau[(x_{1} + v_{1}), (x_{2} + v_{2}), (x_{3} + v_{3})], (\alpha_{1}(x_{4}) + \beta_{1}(v_{4})), (\alpha_{2}(x_{5}) + \beta_{2}(v_{5}))]\cr
&=&\tau[(\alpha_{1}(x_{1}) + \beta_{1}(v_{1})), \tau[(x_{2} + v_{2}), (x_{3} + v_{3}), (x_{4} + v_{4})], (\alpha_{2}(x_{5}) + \alpha_{2}(v_{5}))]\cr
&=&\tau[(\alpha_{1}(x_{1}) + \beta_{1}(v_{1})), (\alpha_{2}(x_{2}) + \beta_{2}(v_{2})), \tau[(x_{3} + v_{3}), (x_{4} + v_{4}), (x_{5} + v_{5})]]. 
\eeqs
After computation, we obtain the Eqs.(\ref{tr1})-(\ref{tr6}). Then ($ \mathcal{L}_{\mu}, \mathcal{M}_{\mu}, \mathcal{R}_{\mu}, \beta_{1}, \beta_{2}, V$) is a quasi trimodule of totally hom-associative ternary algebra $(\A, \mu, \alpha_{1}, \alpha_{2})$ if and only if  $(\mathcal{A} \oplus V, \tau, \alpha_{1} + \beta_{1}, \alpha_{2} + \beta_{2}) $ is a totally hom-associative ternary algebra.

$ \hfill \square $

\begin{example}\label{extha}
 Let $(\mathcal{A}, \mu, \alpha_{1}, \alpha_{2})$ be a totally multiplicative hom-associative ternary algebra.

 Considering the linear maps of the  left $L_{\mu},$ right $R_{\mu}$ and central $M_{\mu}$ multiplication operations defined as:
 \begin{eqnarray}
 L_{\mu}: \A\otimes \A & \longrightarrow & \mathfrak{gl}(\A)  \cr
  (x, y)  & \longmapsto & L_{\mu}(x, y):
  \begin{array}{ccc}
 \A &\longrightarrow & \A \cr 
  z & \longmapsto & \mu(x, y, z), 
   \end{array}
\end{eqnarray}
\begin{eqnarray}
    R_{\mu}: \A\otimes \A & \longrightarrow & \mathfrak{gl}(\A)  \cr
     (x, y)  & \longmapsto & R_{\mu}(x, y):
     \begin{array}{ccc}
    \A &\longrightarrow & \A \cr 
     z & \longmapsto & \mu(z, x, y),
      \end{array}
\end{eqnarray}
      \begin{eqnarray}
    M_{\mu}: \A\otimes \A & \longrightarrow & \mathfrak{gl}(\A)  \cr
     (x, y)  & \longmapsto & M_{\mu}(x, y):
     \begin{array}{ccc}
    \A &\longrightarrow & \A \cr 
     z & \longmapsto & \mu(x, z, y).
      \end{array}
 \end{eqnarray}
Then ($ L_{\mu}, 0, 0, \alpha_{1}, \alpha_{2}, \A$), ($ 0, 0, R_{\mu}, \alpha_{1}, \alpha_{2}, \A$) and ($ L_{\mu}, M_{\mu}, R_{\mu}, \alpha_{1}, \alpha_{2}, \A$) are quasi trimodules of totally hom-associative ternary algebra $(\A, \mu, \alpha_{1}, \alpha_{2})$.
\end{example}

\begin{theorem}\label{Theo. Match. Pair for totally}
 Let $(\A, \mu_{\A}, \alpha_{1}, \alpha_{2})$  and $(\B, \mu_{\B}, \beta_{1}, \beta_{2})$ be two totally hom-associative ternary algebras. Suppose that
 there are linear maps $\mathcal{L}_{\mu_{\A}}: \A\otimes\A\otimes\B\rightarrow \B, 
 \mathcal{R}_{\mu_{\A}}: \B\otimes\A\otimes\A\rightarrow \B, \mathcal{M}_{\mu_{\A}}: 
 \A\otimes\B\otimes\A\rightarrow \B$ and $\mathcal{L}_{\mu_{\B}}: \B\otimes\B\otimes\A\rightarrow 
 \A, \mathcal{R}_{\mu_{\B}}: \A\otimes\B\otimes\B\rightarrow \A, \mathcal{M}_{\mu_{\B}}: 
 \B\otimes\A\otimes\B\rightarrow \A$ such that ($ \mathcal{L}_{\mu_{\A}}, \mathcal{M}_{\mu_{\A}}, 
 \mathcal{R}_{\mu_{\A}}, \beta_{1}, \beta_{2}, \B$) is a quasi trimodule of the totally hom-associative ternary algebra 
 $(\mathcal{A}, \mu_{A}, \alpha_{1}, \alpha_{2})$ and ($ \mathcal{L}_{\mu_{\B}}, \mathcal{M}_{\mu_{\B}}, 
 \mathcal{R}_{\mu_{\B}}, \alpha_{1}, \alpha_{2}, \A$) is a quasi trimodule of the totally hom-associative ternary algebra 
 $(\mathcal{B}, \mu_{B}, \beta_{1}, \beta_{2}),$  satisfying the following conditions:
 \beq \label{tr7}
\mu_{\A}(\mathcal{L}_{\mu_{\B}}(a, b)(x), \alpha_{1}(y), \alpha_{2}(z))&=& \mathcal{L}_{\mu_{\B}}[\beta_{1}(a), \mathcal{R}_{\mu_{\A}}(x, y)(b)]\alpha_{2}(z)\cr
&=& \mathcal{L}_{\mu_{\B}}(\beta_{1}(a), \beta_{2}(b))(\mu_{\A}(x, y, z)),
\eeq
 \beq \label{tr8}
\mu_{\A}(\mathcal{M}_{\mu_{\B}}(a, b)(x), \alpha_{1}(y), \alpha_{2}(z))&=& \mathcal{L}_{\mu_{\B}}[\beta_{1}(a), \mathcal{M}_{\mu_{\A}}(x, y)(b)]\alpha_{2}(z)\cr
&=& \mathcal{M}_{\mu_{\B}}[\beta_{1}(a), \mathcal{R}_{\mu_{\A}}(y, z)(b)]\alpha_{2}(x),
\eeq
 \beq \label{tr9}
\mu_{\A}(\mathcal{R}_{\mu_{\B}}(a, b)(x), \alpha_{1}(y), \alpha_{2}(z))&=& \mu_{\A}[\alpha_{1}(x), \mathcal{L}_{\mu_{\B}}(a, b)(y), \alpha_{2}(z)]\cr
&=& \mathcal{R}_{\mu_{\B}}[\beta_{2}(a), \mathcal{R}_{\mu_{\A}}(y, z)(b)]\alpha_{1}(x),
\eeq
\beq \label{tr10}
\mathcal{L}_{\mu_{\B}}[\mathcal{L}_{\mu_{\A}}(x, y)(a), \beta_{1}(b)]\alpha_{2}(z)&=& \mu_{\A}[\alpha_{1}(x), \mathcal{R}_{\mu_{\B}}(a, b)(y), \alpha_{2}(z)]\cr
&=& \mu_{\A}[\alpha_{1}(x), \alpha_{2}( y), \mathcal{L}_{\mu_{\B}}(a, b)(z)],
\eeq
\beq \label{tr11}
\mathcal{L}_{\mu_{\B}}[\mathcal{M}_{\mu_{\A}}(x, y)(a), \beta_{1}(b)]\alpha_{2}(z)&=& \mu_{\A}[\alpha_{1}(x), \mathcal{M}_{\mu_{\B}}(a, b)(y), \alpha_{2}(z)]\cr
&=& \mathcal{R}_{\mu_{\B}}[\beta_{2}(a), \mathcal{M}_{\mu_{\A}}(y, z)(b)]\alpha_{1}(x),
\eeq
\beq \label{tr12}
\mathcal{L}_{\mu_{\B}}[\mathcal{R}_{\mu_{\A}}(x, y)(a), \beta_{1}(b)]\alpha_{2}(z)&=& \mathcal{L}_{\mu_{\B}}[\beta_{1}(a), \mathcal{L}_{\mu_{\A}}(x, y)(b)]\alpha_{2}(z)\cr
&=& \mathcal{M}_{\mu_{\B}}[\beta_{1}(a), \mathcal{M}_{\mu_{\A}}(y, z)(b)]\alpha_{2}(x),
\eeq
\beq \label{tr13}
\mathcal{M}_{\mu_{\B}}[\mathcal{L}_{\mu_{\A}}(x, y)(a), \beta_{2}(b)]\alpha_{1}(z)&=& \mathcal{R}_{\mu_{\B}}[ \mathcal{M}_{\mu_{\A}}(y, z)(a), \beta_{2}(b)]\alpha_{1}(x)\cr
&=& \mu_{\A}[\alpha_{1}(x), \alpha_{2}(y), \mathcal{M}_{\mu_{\B}}(a, b)(z)],
\eeq
\beq \label{tr14}
\mathcal{M}_{\mu_{\B}}[\mathcal{M}_{\mu_{\A}}(x, y)(a), \beta_{2}(b)]\alpha_{1}(z)&=& \mathcal{R}_{\mu_{\B}}[ \mathcal{R}_{\mu_{\A}}(y, z)(a), \beta_{2}(b)]\alpha_{1}(x)\cr
&=& \mathcal{R}_{\mu_{\B}}[\beta_{2}(a), \mathcal{L}_{\mu_{\A}}(y, z)(b)]\alpha_{1}(x),
\eeq
\beq \label{tr15}
\mathcal{M}_{\mu_{\B}}[\mathcal{R}_{\mu_{\A}}(x, y)(a), \beta_{2}(b)]\alpha_{1}(z)&=& \mathcal{M}_{\mu_{\B}}(\beta_{1}(a), \beta_{2}(b)) (\mu_{\A}(x, y, z)\cr
&=& \mathcal{M}_{\mu_{\B}}[\beta_{1}(a), \mathcal{L}_{\mu_{\A}}(y, z)(b)]\alpha_{2}(x),
\eeq
\beq \label{tr16}
\mathcal{R}_{\mu_{\B}}(\beta_{1}(a), \beta_{2}(b))(\mu_{\A}(x, y, z))&=& \mathcal{R}_{\mu_{\B}}[\mathcal{L}_{\mu_{\A}}(y, z)(a), \beta_{2}(b)]\alpha_{1}(x)\cr
&=& \mu_{\A}[\alpha_{1}(x), \alpha_{2}(y), \mathcal{R}_{\mu_{\B}}(a, b)(z)],
\eeq
\beq \label{tr17}
\mu_{\B}[\mathcal{L}_{\mu_{\A}}(x, y)(a), \beta_{1}(b), \beta_{2}(c)]&=& \mathcal{L}_{\mu_{\A}}[\alpha_{1}(x), \mathcal{R}_{\mu_{\B}}(a, b)(y)]\beta_{2}((c)\cr
&=& \mathcal{L}_{\mu_{\A}}(\alpha_{1}(x), \alpha_{2}(y))(\mu_{\B}(a, b, c)),
\eeq
\beq \label{tr18}
\mu_{\B}[\mathcal{M}_{\mu_{\A}}(x, y)(a), \beta_{1}(b), \beta_{2}(c)]&=& \mathcal{L}_{\mu_{\A}}[\alpha_{1}(x), \mathcal{M}_{\mu_{\B}}(a, b)(y)]\beta_{2}(c)\cr
&=& \mathcal{M}_{\mu_{\A}}[\alpha_{1}(x), \mathcal{R}_{\mu_{\B}}(b, c)(y)]\beta_{2}(a),
\eeq
\beq \label{tr19}
\mu_{\B}[\mathcal{R}_{\mu_{\A}}(x, y)(a), \beta_{1}(b), \beta_{2}(c)]&=& \mu_{\B}[\beta_{1}(a), \mathcal{L}_{\mu_{\A}}(x, y)(b), \beta_{2}(c)]\cr
&=& \mathcal{R}_{\mu_{\A}}[\alpha_{2}(x), \mathcal{R}_{\mu_{\B}}(b, c)(y)]\beta_{1}(a),
\eeq
\beq \label{tr20}
\mathcal{L}_{\mu_{\A}}[\mathcal{L}_{\mu_{\B}}(a, b)(x), \alpha_{1}(y)]\beta_{2}(c)&=& \mu_{\B}[\beta_{1}(a), \mathcal{R}_{\mu_{\A}}(x, y)(b), \beta_{2}(c)]\cr
&=& \mu_{\B}[\beta_{1}(a), \beta_{2}(b), \mathcal{L}_{\mu_{\A}}(x, y)(c)],
\eeq
\beq \label{tr21}
\mathcal{L}_{\mu_{\A}}[\mathcal{M}_{\mu_{\B}}(a, b)(x), \alpha_{1}(y)]\beta_{2}(c)&=& \mu_{\B}[\beta_{1}(a), \mathcal{M}_{\mu_{\A}}(x, y)(b), \beta_{2}(c)]\cr 
&=& \mathcal{R}_{\mu_{\A}}[\alpha_{2}(x), \mathcal{M}_{\mu_{\B}}(b, c)(y)]\beta_{1}(a),
\eeq
\beq \label{tr22}
\mathcal{L}_{\mu_{\A}}[\mathcal{R}_{\mu_{\B}}(a, b)(x), \alpha_{1}(y)]\beta_{2}(c)&=& \mathcal{L}_{\mu_{\A}}[\alpha_{1}(x), \mathcal{L}_{\mu_{\B}}(a, b)(y)]\beta_{2}(c)\cr
&=& \mathcal{M}_{\mu_{\A}}[\alpha_{1}(x), \mathcal{M}_{\mu_{\B}}(b, c)(y)]\beta_{2}(a),
\eeq
\beq \label{tr23}
\mathcal{M}_{\mu_{\A}}[\mathcal{L}_{\mu_{\B}}(a, b)(x), \alpha_{2}(y)]\beta_{1}(c)&=& \mathcal{R}_{\mu_{\A}}[ \mathcal{M}_{\mu_{\B}}(b, c)(x), \alpha_{2}(y)]\beta_{1}(a)\cr
&=& \mu_{\B}[\beta_{1}(a), \beta_{2}(b), \mathcal{M}_{\mu_{\A}}(x, y)(c)],
\eeq
\beq \label{tr24}
\mathcal{M}_{\mu_{\A}}[\mathcal{M}_{\mu_{\B}}(a, b)(x), \alpha_{2}(y)]\beta_{1}(c)&=& \mathcal{R}_{\mu_{\A}}[ \mathcal{R}_{\mu_{\B}}(b, c)(x), \alpha_{2}(y)]\beta_{1}(a)\cr
&=& \mathcal{R}_{\mu_{\A}}[\alpha_{2}(x), \mathcal{L}_{\mu_{\B}}(b, c)(y)]\beta_{1}(a),
\eeq
\beq \label{tr25}
\mathcal{M}_{\mu_{\A}}[\mathcal{R}_{\mu_{\B}}(a, b)(x), \alpha_{2}(y)]\beta_{1}(c)&=& \mathcal{M}_{\mu_{\A}}(\alpha_{1}(x), \alpha_{2}(y))(\mu_{\B}(a, b, c))\cr
&=& \mathcal{M}_{\mu_{\A}}[\alpha_{1}(x), \mathcal{L}_{\mu_{\B}}(b, c)(y)]\beta_{2}(a),
\eeq
\beq \label{tr26}
\mathcal{R}_{\mu_{\A}}(\alpha_{1}(x), \alpha_{2}(y))(\mu_{\B}(a, b, c))&=& \mathcal{R}_{\mu_{\A}}[\mathcal{L}_{\mu_{\B}}(b, c)(x), \alpha_{2}(y)]\beta_{1}(a)\cr
&=& \mu_{\B}[\beta_{1}(a), \beta_{2}(b), \mathcal{R}_{\mu_{\A}}(x, y)(c)],
\eeq
for any $ x, y, z\in \mathcal{A}, a, b, c\in \mathcal{B}$. Then, there is a totally hom-associative ternary algebra structure on the direct sum $\mathcal{A} \oplus \mathcal{B}$ of the underlying vector spaces of $ \mathcal{A} $ and $\mathcal{B}$ given by the product $\tau$ defined {as}
\beqs
\tau[(x + a), (y + b), (z + c)]&=& [\mu_{\mathcal{A}}(x, y, z) + \mathcal{L}_{\mu_{\B}}(a, b)(z) + \mathcal{M}_{\mu_{\B}}(a, c)(y) + \mathcal{R}_{\mu_{\B}}(b, c)(x)] +\cr
&& [\mu_{\mathcal{B}}(a, b, c) + \mathcal{L}_{\mu_{\A}}(x, y)(c) + \mathcal{M}_{\mu_{\A}}(x, z)(b) + \mathcal{R}_{\mu_{\A}}(y, z)(a)]
\eeqs
for any $ x, y, z\in \mathcal{A}, a, b, c\in \mathcal{B}$.
\end{theorem}
\textbf{Proof}
Let $x_{1}, x_{2}, x_{3}, x_{4}, x_{5}\in \mathcal{A}$ and $y_{1}, y_{2}, y_{3}, y_{4}, y_{5}\in \mathcal{B}$. By definition, we have 
\beqs
\tau[(x + a), (y + b), (z + c)]&=& [\mu_{\mathcal{A}}(x, y, z) + \mathcal{L}_{\mu_{\B}}(a, b)(z) + \mathcal{M}_{\mu_{\B}}(a, c)(y) + \mathcal{R}_{\mu_{\B}}(b, c)(x)] +\cr
&& [\mu_{\mathcal{B}}(a, b, c) + \mathcal{L}_{\mu_{\A}}(x, y)(c) + \mathcal{M}_{\mu_{\A}}(x, z)(b) + \mathcal{R}_{\mu_{\A}}(y, z)(a)]
\eeqs 
for any $x, y, z\in \mathcal{A}, a, b, c\in \mathcal{B}$.
Setting the  conditions:
\beqs
&&\tau[\tau[(x_{1} + y_{1}), (x_{2} + y_{2}), (x_{3} + y_{3})], (\alpha_{1}(x_{4}) + \beta_{2}( y_{4})), (\alpha_{2}(x_{5}) + \beta_{2}(y_{5}))]\cr
&=&\tau[(\alpha_{1}(x_{1}) + \beta_{1}(y_{1})), \tau[(x_{2} + y_{2}), (x_{3} + y_{3}), (x_{4} + y_{4})], (\alpha_{2}(x_{5}) + \beta_{2}(y_{5}))]\cr
&=&\tau[(\alpha_{1}(x_{1}) + \beta_{1}(y_{1})), (\alpha_{2}(x_{2}) + \beta_{2}(y_{2})), \tau[(x_{3} + y_{3}), (x_{4} + y_{4}), (x_{5} + y_{5})]], 
\eeqs
we obtain, by direct computation, the Eqs.(\ref{tr7}) - (\ref{tr26}). Then, there is a totally hom-associative ternary algebra structure on the direct sum $\mathcal{A} \oplus \mathcal{B}$ of the underlying vector spaces of $ \mathcal{A} $ and $ \mathcal{B} $ if and only if the Eqs.(\ref{tr7}) - (\ref{tr26}) are satisfied.
  
$ \hfill \square $

Let $ \mathcal{A} \bowtie^{\mathcal{L}_{\mu_{\mathcal{A}}}, \mathcal{M}_{\mu_{\mathcal{A}}}, 
\mathcal{R}_{\mu_{\mathcal{A}}},  \alpha_{1}, \alpha_{2}}_{\mathcal{L}_{\mu_{\mathcal{B}}}, \mathcal{M}_{\mu_{\mathcal{B}}}, \mathcal{R}_{\mu_{\mathcal{B}}}, \beta_{1}, \beta_{2}} \mathcal{B} $ denote this totally hom-associative ternary algebra.
\begin{definition}
Let $ (\mathcal{A}, \mu_{\mathcal{A}}, \alpha_{1}, \alpha_{2})$ and $  (\mathcal{B}, \mu_{\mathcal{B}}, \beta_{1}, \beta_{2}) $ 
be two totally hom-associative ternary algebras. Suppose there are linear maps 
$\mathcal{L}_{\mu_{\mathcal{A}}}, \mathcal{M}_{\mu_{\mathcal{A}}}, \mathcal{R}_{\mu_{\mathcal{A}}}$ 
and $\mathcal{L}_{\mu_{\mathcal{B}}}, \mathcal{M}_{\mu_{\mathcal{B}}}, \mathcal{R}_{\mu_{\mathcal{B}}}$
 such that\\ $(\mathcal{L}_{\mu_{\mathcal{A}}}, \mathcal{M}_{\mu_{\mathcal{A}}}, \mathcal{R}_{\mu_{\mathcal{A}}}, \beta_{1}, \beta_{2}, \mathcal{B})$ is a quasi trimodule of $(\mathcal{A}, \mu_{\mathcal{A}}, \alpha_{1}, \alpha_{2}),$ 
and $(\mathcal{L}_{\mu_{\mathcal{B}}}, \mathcal{M}_{\mu_{\mathcal{B}}}, \mathcal{R}_{\mu_{\mathcal{B}}}, \alpha_{1}, \alpha_{2}, \mathcal{A})$ is a quasi trimodule of $( \mathcal{B}, \beta_{1}, \beta_{2})$. 
If Eqs.(\ref{tr7})-(\ref{tr26}) are satisfied, then 
\beqs
(\mathcal{A}, \mathcal{B},\mathcal{L}_{\mu_{\mathcal{A}}}, \mathcal{M}_{\mu_{\mathcal{A}}}, \mathcal{R}_{\mu_{\mathcal{A}}}, \alpha_{1}, \alpha_{2}, \mathcal{L}_{\mu_{\mathcal{B}}}, \mathcal{M}_{\mu_{\mathcal{B}}}, \mathcal{R}_{\mu_{\mathcal{B}}}, \beta_{1}, \beta_{2})
\eeqs
 is called a matched pair of totally hom-associative ternary algebras.   
\end{definition}

\begin{definition}
Let $ (\mathcal{A}, \mu_{\mathcal{A}}, \alpha_{1}, \alpha_{2}) $ and $(\mathcal{B}, \mu_{\mathcal{B}}, \beta_{1}, \beta_{2}) $ 
be two totally hom-associative ternary algebras. Suppose there are linear maps 
$\mathcal{L}_{\mu_{\mathcal{A}}}, \mathcal{M}_{\mu_{\mathcal{A}}}, \mathcal{R}_{\mu_{\mathcal{A}}}$ 
and $\mathcal{L}_{\mu_{\mathcal{B}}}, \mathcal{M}_{\mu_{\mathcal{B}}}, \mathcal{R}_{\mu_{\mathcal{B}}}$
 such that\\$(\mathcal{L}_{\mu_{\mathcal{A}}}, \mathcal{M}_{\mu_{\mathcal{A}}}, 
 \mathcal{R}_{\mu_{\mathcal{A}}}, \beta_{1}, \beta_{2}, \mathcal{B})$ is a trimodule of $(\mathcal{A}, 
 \mu_{\mathcal{A}}, \alpha_{1}, \alpha_{2}),$ 
and $(\mathcal{L}_{\mu_{\mathcal{B}}}, \mathcal{M}_{\mu_{\mathcal{B}}}, \mathcal{R}_{\mu_{\mathcal{B}}}, 
\alpha_{1}, \alpha_{2}, \mathcal{A})$ is a trimodule of $ (\mathcal{B}, \mu_{\mathcal{B}}, \beta_{1}, 
\beta_{2})$. 
  Then $(\mathcal{A}, \mathcal{B},\mathcal{L}_{\mu_{\mathcal{A}}}, \mathcal{M}_{\mu_{\mathcal{A}}}, \mathcal{R}_{\mu_{\mathcal{A}}}, \alpha_{1}, \alpha_{2}, \mathcal{L}_{\mu_{\mathcal{B}}}, \mathcal{M}_{\mu_{\mathcal{B}}}, \mathcal{R}_{\mu_{\mathcal{B}}}, \beta_{1}, \beta_{2})$ 
 is a matched pair of totally hom-associative ternary algebras if Eqs.(\ref{tr7})-(\ref{tr26}) are satisfied and  the following conditions hold:
 \beq \label{tr27}
 \mathcal{M}_{\mu_{A}}(\alpha_{1}(a), \alpha_{2}(z))(\mathcal{M}_{\mu_{A}}(\alpha_{1}(b), \alpha_{2}(y))(\mathcal{M}_{\mu_{A}}(\alpha_{1}(c), \alpha_{2}(x))\beta_{i}(v'))=\cr
  \mathcal{M}_{\mu_{A}}(\mu_{A}(\alpha_{1}(a), \alpha_{1}(b), \alpha_{1}(c)), \mu_{A}(\alpha_{2}(x), \alpha_{2}(y), \alpha_{2}(z))\beta_{i}(v'),
 \eeq   
 \beq \label{tr28}
 \mathcal{M}_{\mu_{B}}(\beta_{1}(a'), \beta_{2}(z'))(\mathcal{M}_{\mu_{B}}(\beta_{1}(b'), \beta_{2}(y'))(\mathcal{M}_{\mu_{B}}(\beta_{1}(c'), \beta_{2}(x'))\alpha_{i}(v))=\cr
  \mathcal{M}_{\mu_{B}}(\mu_{B}(\beta_{1}(a'), \beta_{1}(b'), \beta_{1}(c')), \mu_{B}(\beta_{2}(x'), \beta_{2}(y'), \beta_{2}(z')))\alpha_{i}(v),
 \eeq
 \beq\label{tr29}
\beta_{1}(\mathcal{L}_{\mu_{\mathcal{A}}}(a, b)v')= \mathcal{L}_{\mu_{\mathcal{A}}}(\alpha_{1}(a), \alpha_{2}(b))\beta_{1}(v'),\ \beta_{2}(\mathcal{L}_{\mu_{\mathcal{A}}}(a, b)v')= \mathcal{L}_{\mu_{\mathcal{A}}}(\alpha_{1}(a), \alpha_{2}(b))\beta_{2}(v')
\eeq
\beq\label{tr30}
\beta_{1}(\mathcal{M}_{\mu_{\mathcal{A}}}(a, b)v')= \mathcal{M}_{\mu_{\mathcal{A}}}(\alpha_{1}(a), \alpha_{2}(b))\beta_{1}(v'),\
\beta_{2}(\mathcal{M}_{\mu_{\mathcal{A}}}(a, b)v')= \mathcal{M}_{\mu_{\mathcal{A}}}(\alpha_{1}(a), \alpha_{2}(b))\beta_{2}(v')
\eeq
\beq\label{tr31}
\beta_{1}(\mathcal{R}_{\mu_{\mathcal{A}}}(a, b)v')= \mathcal{R}_{\mu_{\mathcal{A}}}(\alpha_{1}(a), \alpha_{2}(b))\beta_{1}(v'), 
\beta_{2}(\mathcal{R}_{\mu_{\mathcal{A}}}(a, b)v')= \mathcal{R}_{\mu_{\mathcal{A}}}(\alpha_{1}(a), \alpha_{2}(b))\beta_{2}(v')
\eeq
 \beq\label{tr32}
\alpha_{1}(\mathcal{L}_{\mu_{\mathcal{B}}}(a', b')v)= \mathcal{L}_{\mu_{\mathcal{B}}}(\beta_{1}(a'), \beta_{2}(b'))\alpha_{1}(v),\ \alpha_{2}(\mathcal{L}_{\mu_{\mathcal{B}}}(a', b')v)= \mathcal{L}_{\mu_{\mathcal{B}}}(\beta_{1}(a'), \beta_{2}(b'))\alpha_{2}(v)
\eeq
 \beq\label{tr33}
\alpha_{1}(\mathcal{M}_{\mu_{\mathcal{B}}}(a', b')v)= \mathcal{M}_{\mu_{\mathcal{B}}}(\beta_{1}(a'), \beta_{2}(b'))\alpha_{1}(v),\ \alpha_{2}(\mathcal{M}_{\mu_{\mathcal{B}}}(a', b')v)= \mathcal{M}_{\mu_{\mathcal{B}}}(\beta_{1}(a'), \beta_{2}(b'))\alpha_{2}(v)
\eeq
 \beq\label{tr34}
\alpha_{1}(\mathcal{R}_{\mu_{\mathcal{B}}}(a', b')v)= \mathcal{R}_{\mu_{\mathcal{B}}}(\beta_{1}(a'), \beta_{2}(b'))\alpha_{1}(v),\ \alpha_{2}(\mathcal{R}_{\mu_{\mathcal{B}}}(a', b')v)= \mathcal{R}_{\mu_{\mathcal{B}}}(\beta_{1}(a'), \beta_{2}(b'))\alpha_{2}(v)
\eeq
 for any $ x, y, z, a, b, c, v\in \mathcal{A}, x', y', z', a', b', c', v'\in \mathcal{B}$ and $i= 1, 2$.
\end{definition}

\subsection{Trimodules and matched pairs of partially hom-associative ternary  algebras}
\begin{definition}\label{trimodule of PHATA}
A trimodule structure over partially hom-associative ternary algebra $(\A, \mu, \alpha_{1}, \alpha_{2})$ on a bihom-module $(V, \beta_{1}, \beta_{2})$ is defined by the following three linear multiplication mappings: 
\beqs
&&\mathcal{L}_{\mu}: \A\otimes\A\otimes V\rightarrow V, \cr
&& \mathcal{R}_{\mu}: V\otimes\A\otimes\A\rightarrow V, \cr
&& \mathcal{M}_{\mu}: \A\otimes V\otimes\A\rightarrow V
\eeqs 
 satisfying the following compatibility conditions
\beq \label{pr1}
\mathcal{L}_{\mu}(\alpha_{1}(a), \alpha_{2}(b))(\mathcal{L}_{\mu}(c, d)(v)) + \mathcal{L}_{\mu}(\mu(a, b, c), \alpha_{1}(d))\beta_{2}(v) + \mathcal{L}_{\mu}(\alpha_{1}(a), \mu(b, c, d))\beta_{2}( v)= 0,
\eeq
\beq \label{pr2}
\mathcal{R}_{\mu}(\alpha_{1}(c), \alpha_{2}(d))(\mathcal{R}_{\mu}(a, b)(v)) + \mathcal{R}_{\mu}(\alpha_{2}(a), \mu(b, c, d))\beta_{1}(v) + \mathcal{R}_{\mu}(\mu(a, b, c), \alpha_{2}(d))\beta_{1}(v)= 0,
\eeq 
\beq \label{pr4}
\mathcal{M}_{\mu}(\alpha_{1}(a), \alpha_{2}(d))(\mathcal{L}{\mu}(b, c)(v)) + \mathcal{L}_{\mu}(\alpha_{1}(a), \alpha_{2}(b)) (\mathcal{M}_{\mu}(c, d)(v))\cr + \M_{\mu}(\mu(a, b, c), \alpha_{2}(d))\beta_{1}(v)= 0,
\eeq
\beq \label{pr5}
\mathcal{M}_{\mu}(\alpha_{1}(a), \alpha_{2}(d))(\mathcal{R}_{\mu}(b, c)(v)) + \mathcal{R}_{\mu}(\alpha_{1}(c), \alpha_{2}(d))(\mathcal{M}_{\mu}(a, b)(v))\cr + \mathcal{M}_{\mu}(\alpha_{1}(a), \mu(b, c, d))\beta_{2}(v)=0,
\eeq
\beq \label{pr6}
\mathcal{R}_{\mu}(\alpha_{1}(c), \alpha_{2}(d))(\mathcal{L}_{\mu}(a, b) (v))+ \mathcal{L}_{\mu}(\alpha_{1}(a), \alpha_{2}(b))(\mathcal{R}_{\mu}(c, d)(v))\cr
+ \mathcal{M}_{\mu}(\alpha_{1}(a), \alpha_{2}(d))(\mathcal{M}_{\mu}(b, c)(v))= 0,
\eeq
\beq \label{pr3}
\mathcal{M}_{\mu}(\alpha_{1}(a), \alpha_{2}(z))(\mathcal{M}_{\mu}(\alpha_{1}(b), \alpha_{2}(y))(\mathcal{M}_{\mu}(\alpha_{1}(c), \alpha_{2}(x))\beta_{1}(v))=\cr
 \mathcal{M}_{\mu}(\mu(\alpha_{1}(a), \alpha_{1}(b), \alpha_{1}(c)), \mu(\alpha_{2}(x), \alpha_{2}(y), \alpha_{2}(z)))\beta_{1}(v),
\eeq
\beq \label{pr3,2}
\mathcal{M}_{\mu}(\alpha_{1}(a), \alpha_{2}(z))(\mathcal{M}_{\mu}(\alpha_{1}(b), \alpha_{2}(y))(\mathcal{M}_{\mu}(\alpha_{1}(c), \alpha_{2}(x))\beta_{2}(v))=\cr
 \mathcal{M}_{\mu}(\mu(\alpha_{1}(a), \alpha_{1}(b), \alpha_{1}(c)), \mu(\alpha_{2}(x), \alpha_{2}(y), \alpha_{2}(z)))\beta_{2}(v),
\eeq
\beq\label{pr7}
\beta_{1}(\mathcal{L}_{\mu}(a, b)v)= \mathcal{L}_{\mu}(\alpha_{1}(a), \alpha_{2}(b))\beta_{1}(v),\ \beta_{2}(\mathcal{L}_{\mu}(a, b)v)= \mathcal{L}_{\mu}(\alpha_{1}(a), \alpha_{2}(b))\beta_{2}(v)
\eeq
\beq\label{pr8}
\beta_{1}(\mathcal{M}_{\mu}(a, b)v)= \mathcal{M}_{\mu}(\alpha_{1}(a), \alpha_{2}(b))\beta_{1}(v),\
\beta_{2}(\mathcal{M}_{\mu}(a, b)v)= \mathcal{M}_{\mu}(\alpha_{1}(a), \alpha_{2}(b))\beta_{2}(v)
\eeq
\beq\label{pr9}
\beta_{1}(\mathcal{R}_{\mu}(a, b)v)= \mathcal{R}_{\mu}(\alpha_{1}(a), \alpha_{2}(b))\beta_{1}(v), 
\beta_{2}(\mathcal{R}_{\mu}(a, b)v)= \mathcal{R}_{\mu}(\alpha_{1}(a), \alpha_{2}(b))\beta_{2}(v)
\eeq
$ \forall a, b, c, d, x, y, z\in \A, v\in V$.
\end{definition}
\begin{remark}
In the case where the Eqs.(\ref{pr1})-(\ref{pr6}) are only satisfied, we refer to the name {of} quasi  trimodule structure instead of simply trimodule structure.
\end{remark}

\begin{proposition}
($ \mathcal{L}_{\mu}, \mathcal{M}_{\mu}, \mathcal{R}_{\mu}, \beta_{1}, \beta_{2}, V$) is a quasi trimodule of partially hom-associative ternary algebra $(\A, \mu, \alpha_{1}, \alpha_{2})$ if and only if the direct sum $(\mathcal{A} \oplus V, \tau, \alpha_{1} + \beta_{1}, \alpha_{2} + \beta_{2}) $ of the underlying vector spaces of $ \mathcal{A} $ and $V$ is turned into a partially hom-associative ternary algebra {product} $\tau$ given by
\beqs
\tau[(x + a), (y + b), (z + c)]= \mu(x, y, z) + \mathcal{L}_{\mu}(x, y)(c) + \mathcal{M}_{\mu}(x, z)(b) + \mathcal{R}_{\mu}(y, z)(a),
\eeqs
for all $ x, y, z \in \mathcal{A}, a, b, c \in V$. We denote it by $ \mathcal{A} \ltimes_{\mathcal{L}_{\mu}, \mathcal{M}_{\mu}, \mathcal{R}_{\mu}, \beta_{1}, \beta_{2}} V$.
\end{proposition}
\textbf{Proof}
Let $v_{1}, v_{2}, v_{3}, v_{4}, v_{5}\in V$ and $x_{1}, x_{2}, x_{3}, x_{4}, x_{5}\in \mathcal{A}$.
Set
\beqs
\tau[\tau[(x_{1} + v_{1}), (x_{2} + v_{2}), (x_{3} + v_{3})], (\alpha_{1}(x_{4}) + \beta_{1}(v_{4})), (\alpha_{2}(x_{5}) + \beta_{2}(v_{5}))]\cr
+\tau[(\alpha_{1}(x_{1}) + \beta_{1}(v_{1})), \tau[(x_{2} + v_{2}), (x_{3} + v_{3}), (x_{4} + v_{4})], (\alpha_{2}(x_{5}) + \alpha_{2}(v_{5}))]\cr
+\tau[(\alpha_{1}(x_{1}) + \beta_{1}(v_{1})), (\alpha_{2}(x_{2}) + \beta_{2}(v_{2})), \tau[(x_{3} + v_{3}), (x_{4} + v_{4}), (x_{5} + v_{5})]]\cr
=0. 
\eeqs
After computation, we obtain the Eqs.(\ref{tr1})-(\ref{tr6}). Then ($ \mathcal{L}_{\mu}, \mathcal{M}_{\mu}, \mathcal{R}_{\mu}, \beta_{1}, \beta_{2}, V$) is a quasi trimodule of partially hom-associative ternary algebra $(\A, \mu, \alpha_{1}, \alpha_{2})$ if and only if  $(\mathcal{A} \oplus V, \tau, \alpha_{1} + \beta_{1}, \alpha_{2} + \beta_{2}) $ is a partially hom-associative ternary algebra.

$ \hfill \square $

\begin{example}
 Let $(\mathcal{A}, \mu, \alpha_{1}, \alpha_{2})$ be a partially multiplicative hom-associative ternary algebra.
 Consider the linear maps of the  left $L_{\mu},$ right $R_{\mu}$ and central $M_{\mu}$ {multiplications}. Then ($ L_{\mu}, 0, 0, \alpha_{1}, \alpha_{2}, \A$), ($ 0, 0, R_{\mu}, \alpha_{1}, \alpha_{2}, \A$) and ($ L_{\mu}, M_{\mu}, R_{\mu}, \alpha_{1}, \alpha_{2}, \A$) are quasi trimodules of partially hom-associative ternary algebra $(\A, \mu, \alpha_{1}, \alpha_{2})$.
\end{example}

\begin{theorem}\label{Theo. Match. Pair for partially}
 Let $(\A, \mu_{\A}, \alpha_{1}, \alpha_{2})$  and $(\B, \mu_{\B}, \beta_{1}, \beta_{2})$ be two partially hom-associative ternary algebras. Suppose that
 there are linear maps $\mathcal{L}_{\mu_{\A}}: \A\otimes\A\otimes\B\rightarrow \B, 
 \mathcal{R}_{\mu_{\A}}: \B\otimes\A\otimes\A\rightarrow \B, \mathcal{M}_{\mu_{\A}}: 
 \A\otimes\B\otimes\A\rightarrow \B,$ and $\mathcal{L}_{\mu_{\B}}: \B\otimes\B\otimes\A\rightarrow 
 \A, \mathcal{R}_{\mu_{\B}}: \A\otimes\B\otimes\B\rightarrow \A, \mathcal{M}_{\mu_{\B}}: 
 \B\otimes\A\otimes\B\rightarrow \A$ such that ($ \mathcal{L}_{\mu_{\A}}, \mathcal{M}_{\mu_{\A}}, 
 \mathcal{R}_{\mu_{\A}}, \beta_{1}, \beta_{2}, \B$) is a quasi trimodule of the partially hom-associative ternary algebra 
 $(\mathcal{A}, \mu_{A}, \alpha_{1}, \alpha_{2}),$ and ($ \mathcal{L}_{\mu_{\B}}, \mathcal{M}_{\mu_{\B}}, 
 \mathcal{R}_{\mu_{\B}}, \alpha_{1}, \alpha_{2}, \A$) is a quasi trimodule of the partially hom-associative ternary algebra 
 $(\mathcal{B}, \mu_{B}, \beta_{1}, \beta_{2}),$ {satisfying} the following conditions:
 \beq \label{pr7}
\mu_{\A}(\mathcal{L}_{\mu_{\B}}(a, b)(x), \alpha_{1}(y), \alpha_{2}(z))+ \mathcal{L}_{\mu_{\B}}[\beta_{1}(a), \mathcal{R}_{\mu_{\A}}(x, y)(b)]\alpha_{2}(z)\cr
+ \mathcal{L}_{\mu_{\B}}(\beta_{1}(a), \beta_{2}(b))(\mu_{\A}(x, y, z))=0,
\eeq
 \beq \label{pr8}
\mu_{\A}(\mathcal{M}_{\mu_{\B}}(a, b)(x), \alpha_{1}(y), \alpha_{2}(z))+ \mathcal{L}_{\mu_{\B}}[\beta_{1}(a), \mathcal{M}_{\mu_{\A}}(x, y)(b)]\alpha_{2}(z)\cr
+ \mathcal{M}_{\mu_{\B}}[\beta_{1}(a), \mathcal{R}_{\mu_{\A}}(y, z)(b)]\alpha_{2}(x)=0,
\eeq
 \beq \label{pr9}
\mu_{\A}(\mathcal{R}_{\mu_{\B}}(a, b)(x), \alpha_{1}(y), \alpha_{2}(z))+ \mu_{\A}[\alpha_{1}(x), \mathcal{L}_{\mu_{\B}}(a, b)(y), \alpha_{2}(z)]\cr
+\mathcal{R}_{\mu_{\B}}[\beta_{2}(a), \mathcal{R}_{\mu_{\A}}(y, z)(b)]\alpha_{1}(x)=0,
\eeq
\beq \label{pr10}
\mathcal{L}_{\mu_{\B}}[\mathcal{L}_{\mu_{\A}}(x, y)(a), \beta_{1}(b)]\alpha_{2}(z)+ \mu_{\A}[\alpha_{1}(x), \mathcal{R}_{\mu_{\B}}(a, b)(y), \alpha_{2}(z)]\cr
+ \mu_{\A}[\alpha_{1}(x), \alpha_{2}( y), \mathcal{L}_{\mu_{\B}}(a, b)(z)]=0,
\eeq
\beq \label{pr11}
\mathcal{L}_{\mu_{\B}}[\mathcal{M}_{\mu_{\A}}(x, y)(a), \beta_{1}(b)]\alpha_{2}(z)+ \mu_{\A}[\alpha_{1}(x), \mathcal{M}_{\mu_{\B}}(a, b)(y), \alpha_{2}(z)]\cr
+ \mathcal{R}_{\mu_{\B}}[\beta_{2}(a), \mathcal{M}_{\mu_{\A}}(y, z)(b)]\alpha_{1}(x)=0,
\eeq
\beq \label{pr12}
\mathcal{L}_{\mu_{\B}}[\mathcal{R}_{\mu_{\A}}(x, y)(a), \beta_{1}(b)]\alpha_{2}(z)+ \mathcal{L}_{\mu_{\B}}[\beta_{1}(a), \mathcal{L}_{\mu_{\A}}(x, y)(b)]\alpha_{2}(z)\cr
+ \mathcal{M}_{\mu_{\B}}[\beta_{1}(a), \mathcal{M}_{\mu_{\A}}(y, z)(b)]\alpha_{2}(x)=0,
\eeq
\beq \label{pr13}
\mathcal{M}_{\mu_{\B}}[\mathcal{L}_{\mu_{\A}}(x, y)(a), \beta_{2}(b)]\alpha_{1}(z)+ \mathcal{R}_{\mu_{\B}}[ \mathcal{M}_{\mu_{\A}}(y, z)(a), \beta_{2}(b)]\alpha_{1}(x)\cr
+ \mu_{\A}[\alpha_{1}(x), \alpha_{2}(y), \mathcal{M}_{\mu_{\B}}(a, b)(z)]=0,
\eeq
\beq \label{pr14}
\mathcal{M}_{\mu_{\B}}[\mathcal{M}_{\mu_{\A}}(x, y)(a), \beta_{2}(b)]\alpha_{1}(z)+ \mathcal{R}_{\mu_{\B}}[ \mathcal{R}_{\mu_{\A}}(y, z)(a), \beta_{2}(b)]\alpha_{1}(x)\cr
+ \mathcal{R}_{\mu_{\B}}[\beta_{2}(a), \mathcal{L}_{\mu_{\A}}(y, z)(b)]\alpha_{1}(x)=0,
\eeq
\beq \label{pr15}
\mathcal{M}_{\mu_{\B}}[\mathcal{R}_{\mu_{\A}}(x, y)(a), \beta_{2}(b)]\alpha_{1}(z)+ \mathcal{M}_{\mu_{\B}}(\beta_{1}(a), \beta_{2}(b)) (\mu_{\A}(x, y, z)\cr
+ \mathcal{M}_{\mu_{\B}}[\beta_{1}(a), \mathcal{L}_{\mu_{\A}}(y, z)(b)]\alpha_{2}(x)=0,
\eeq
\beq \label{pr16}
\mathcal{R}_{\mu_{\B}}(\beta_{1}(a), \beta_{2}(b))(\mu_{\A}(x, y, z))+ \mathcal{R}_{\mu_{\B}}[\mathcal{L}_{\mu_{\A}}(y, z)(a), \beta_{2}(b)]\alpha_{1}(x)\cr
+ \mu_{\A}[\alpha_{1}(x), \alpha_{2}(y), \mathcal{R}_{\mu_{\B}}(a, b)(z)]=0,
\eeq
\beq \label{pr17}
\mu_{\B}[\mathcal{L}_{\mu_{\A}}(x, y)(a), \beta_{1}(b), \beta_{2}(c)]+ \mathcal{L}_{\mu_{\A}}[\alpha_{1}(x), \mathcal{R}_{\mu_{\B}}(a, b)(y)]\beta_{2}((c)\cr
+ \mathcal{L}_{\mu_{\A}}(\alpha_{1}(x), \alpha_{2}(y))(\mu_{\B}(a, b, c))=0,
\eeq
\beq \label{pr18}
\mu_{\B}[\mathcal{M}_{\mu_{\A}}(x, y)(a), \beta_{1}(b), \beta_{2}(c)]+ \mathcal{L}_{\mu_{\A}}[\alpha_{1}(x), \mathcal{M}_{\mu_{\B}}(a, b)(y)]\beta_{2}(c)\cr
+ \mathcal{M}_{\mu_{\A}}[\alpha_{1}(x), \mathcal{R}_{\mu_{\B}}(b, c)(y)]\beta_{2}(a)=0,
\eeq
\beq \label{pr19}
\mu_{\B}[\mathcal{R}_{\mu_{\A}}(x, y)(a), \beta_{1}(b), \beta_{2}(c)]+ \mu_{\B}[\beta_{1}(a), \mathcal{L}_{\mu_{\A}}(x, y)(b), \beta_{2}(c)]\cr
+ \mathcal{R}_{\mu_{\A}}[\alpha_{2}(x), \mathcal{R}_{\mu_{\B}}(b, c)(y)]\beta_{1}(a)=0,
\eeq
\beq \label{pr20}
\mathcal{L}_{\mu_{\A}}[\mathcal{L}_{\mu_{\B}}(a, b)(x), \alpha_{1}(y)]\beta_{2}(c)+ \mu_{\B}[\beta_{1}(a), \mathcal{R}_{\mu_{\A}}(x, y)(b), \beta_{2}(c)]\cr
+ \mu_{\B}[\beta_{1}(a), \beta_{2}(b), \mathcal{L}_{\mu_{\A}}(x, y)(c)]=0,
\eeq
\beq \label{pr21}
\mathcal{L}_{\mu_{\A}}[\mathcal{M}_{\mu_{\B}}(a, b)(x), \alpha_{1}(y)]\beta_{2}(c)+ \mu_{\B}[\beta_{1}(a), \mathcal{M}_{\mu_{\A}}(x, y)(b), \beta_{2}(c)]\cr 
+ \mathcal{R}_{\mu_{\A}}[\alpha_{2}(x), \mathcal{M}_{\mu_{\B}}(b, c)(y)]\beta_{1}(a)=0,
\eeq
\beq \label{pr22}
\mathcal{L}_{\mu_{\A}}[\mathcal{R}_{\mu_{\B}}(a, b)(x), \alpha_{1}(y)]\beta_{2}(c)+ \mathcal{L}_{\mu_{\A}}[\alpha_{1}(x), \mathcal{L}_{\mu_{\B}}(a, b)(y)]\beta_{2}(c)\cr
+ \mathcal{M}_{\mu_{\A}}[\alpha_{1}(x), \mathcal{M}_{\mu_{\B}}(b, c)(y)]\beta_{2}(a)=0,
\eeq
\beq \label{pr23}
\mathcal{M}_{\mu_{\A}}[\mathcal{L}_{\mu_{\B}}(a, b)(x), \alpha_{2}(y)]\beta_{1}(c)+ \mathcal{R}_{\mu_{\A}}[ \mathcal{M}_{\mu_{\B}}(b, c)(x), \alpha_{2}(y)]\beta_{1}(a)\cr
+ \mu_{\B}[\beta_{1}(a), \beta_{2}(b), \mathcal{M}_{\mu_{\A}}(x, y)(c)]=0,
\eeq
\beq \label{pr24}
\mathcal{M}_{\mu_{\A}}[\mathcal{M}_{\mu_{\B}}(a, b)(x), \alpha_{2}(y)]\beta_{1}(c)+ \mathcal{R}_{\mu_{\A}}[ \mathcal{R}_{\mu_{\B}}(b, c)(x), \alpha_{2}(y)]\beta_{1}(a)\cr
+ \mathcal{R}_{\mu_{\A}}[\alpha_{2}(x), \mathcal{L}_{\mu_{\B}}(b, c)(y)]\beta_{1}(a)=0,
\eeq
\beq \label{pr25}
\mathcal{M}_{\mu_{\A}}[\mathcal{R}_{\mu_{\B}}(a, b)(x), \alpha_{2}(y)]\beta_{1}(c)+ \mathcal{M}_{\mu_{\A}}(\alpha_{1}(x), \alpha_{2}(y))(\mu_{\B}(a, b, c))\cr
+ \mathcal{M}_{\mu_{\A}}[\alpha_{1}(x), \mathcal{L}_{\mu_{\B}}(b, c)(y)]\beta_{2}(a)=0,
\eeq
\beq \label{pr26}
\mathcal{R}_{\mu_{\A}}(\alpha_{1}(x), \alpha_{2}(y))(\mu_{\B}(a, b, c))+ \mathcal{R}_{\mu_{\A}}[\mathcal{L}_{\mu_{\B}}(b, c)(x), \alpha_{2}(y)]\beta_{1}(a)\cr
+ \mu_{\B}[\beta_{1}(a), \beta_{2}(b), \mathcal{R}_{\mu_{\A}}(x, y)(c)]=0,
\eeq
for any $ x, y, z\in \mathcal{A}, a, b, c\in \mathcal{B}$. Then, there is a partially hom-associative ternary algebra structure on the direct sum $\mathcal{A} \oplus \mathcal{B}$ of the underlying vector spaces of $ \mathcal{A} $ and $ \mathcal{B} $ given by the product $\tau$ defined by 
\beqs
\tau[(x + a), (y + b), (z + c)]&=& [\mu_{\mathcal{A}}(x, y, z) + \mathcal{L}_{\mu_{\B}}(a, b)(z) + \mathcal{M}_{\mu_{\B}}(a, c)(y) + \mathcal{R}_{\mu_{\B}}(b, c)(x)] +\cr
&& [\mu_{\mathcal{B}}(a, b, c) + \mathcal{L}_{\mu_{\A}}(x, y)(c) + \mathcal{M}_{\mu_{\A}}(x, z)(b) + \mathcal{R}_{\mu_{\A}}(y, z)(a)]
\eeqs
for any $ x, y, z\in \mathcal{A}, a, b, c\in \mathcal{B}$. 
\end{theorem}
\textbf{Proof}
Let $x_{1}, x_{2}, x_{3}, x_{4}, x_{5}\in \mathcal{A}$ and $y_{1}, y_{2}, y_{3}, y_{4}, y_{5}\in \mathcal{B}$. By definition, we have 
\beqs
\tau[(x + a), (y + b), (z + c)]&=& [\mu_{\mathcal{A}}(x, y, z) + \mathcal{L}_{\mu_{\B}}(a, b)(z) + \mathcal{M}_{\mu_{\B}}(a, c)(y) + \mathcal{R}_{\mu_{\B}}(b, c)(x)] +\cr
&& [\mu_{\mathcal{B}}(a, b, c) + \mathcal{L}_{\mu_{\A}}(x, y)(c) + \mathcal{M}_{\mu_{\A}}(x, z)(b) + \mathcal{R}_{\mu_{\A}}(y, z)(a)]
\eeqs 
for any $x, y, z\in \mathcal{A}, a, b, c\in \mathcal{B}$.
Setting the strong condition 
\beqs
\tau[\tau[(x_{1} + y_{1}), (x_{2} + y_{2}), (x_{3} + y_{3})], (\alpha_{1}(x_{4}) + \beta_{2}( y_{4})), (\alpha_{2}(x_{5}) + \beta_{2}(y_{5}))]\cr
+\tau[(\alpha_{1}(x_{1}) + \beta_{1}(y_{1})), \tau[(x_{2} + y_{2}), (x_{3} + y_{3}), (x_{4} + y_{4})], (\alpha_{2}(x_{5}) + \beta_{2}(y_{5}))]\cr
+\tau[(\alpha_{1}(x_{1}) + \beta_{1}(y_{1})), (\alpha_{2}(x_{2}) + \beta_{2}(y_{2})), \tau[(x_{3} + y_{3}), (x_{4} + y_{4}), (x_{5} + y_{5})]]\cr
=0, 
\eeqs
we obtain, by dierct computation, the Eqs.(\ref{pr7}) - (\ref{pr26}). Then, there is a partially hom-associative ternary algebra structure on the direct sum $\mathcal{A} \oplus \mathcal{B}$ of the underlying vector spaces of $ \mathcal{A} $ and $ \mathcal{B} $ if and only if the Eqs.(\ref{pr7}) - (\ref{pr26}) are satisfied.
  
$ \hfill \square $

Let $ \mathcal{A} \bowtie^{\mathcal{L}_{\mu_{\mathcal{A}}}, \mathcal{M}_{\mu_{\mathcal{A}}}, 
\mathcal{R}_{\mu_{\mathcal{A}}},  \alpha_{1}, \alpha_{2}}_{\mathcal{L}_{\mu_{\mathcal{B}}}, \mathcal{M}_{\mu_{\mathcal{B}}}, \mathcal{R}_{\mu_{\mathcal{B}}}, \beta_{1}, \beta_{2}} \mathcal{B} $ denote this partially hom-associative ternary algebra.

\begin{definition}
Let $ (\mathcal{A}, \mu_{\mathcal{A}}, \alpha_{1}, \alpha_{2}) $ and $  (\mathcal{B}, \mu_{\mathcal{B}}, \beta_{1}, \beta_{2}) $ 
be two partially hom-associative ternary algebras. Suppose that there are linear maps 
$\mathcal{L}_{\mu_{\mathcal{A}}}, \mathcal{M}_{\mu_{\mathcal{A}}}, \mathcal{R}_{\mu_{\mathcal{A}}}$ 
and $\mathcal{L}_{\mu_{\mathcal{B}}}, \mathcal{M}_{\mu_{\mathcal{B}}}, \mathcal{R}_{\mu_{\mathcal{B}}}$
 such that $(\mathcal{L}_{\mu_{\mathcal{A}}}, \mathcal{M}_{\mu_{\mathcal{A}}}, \mathcal{R}_{\mu_{\mathcal{A}}}, \beta_{1}, \beta_{2}, \mathcal{B})$ is a quasi trimodule of $(\mathcal{A}, , \alpha_{1}, \alpha_{2}),$ 
and $(\mathcal{L}_{\mu_{\mathcal{B}}}, \mathcal{M}_{\mu_{\mathcal{B}}}, \mathcal{R}_{\mu_{\mathcal{B}}}, \alpha_{1}, \alpha_{2} \mathcal{A})$ is a quasi trimodule of $( \mathcal{B}, \beta_{1}, \beta_{2})$. 
If Eqs.(\ref{pr7})-(\ref{pr26}) are satisfied, then 
\beqs
(\mathcal{A}, \mathcal{B},\mathcal{L}_{\mu_{\mathcal{A}}}, \mathcal{M}_{\mu_{\mathcal{A}}}, \mathcal{R}_{\mu_{\mathcal{A}}}, \alpha_{1}, \alpha_{2}, \mathcal{L}_{\mu_{\mathcal{B}}}, \mathcal{M}_{\mu_{\mathcal{B}}}, \mathcal{R}_{\mu_{\mathcal{B}}}, \beta_{1}, \beta_{2})
\eeqs
 is called a matched pair of partially hom-associative ternary algebras.   
\end{definition}

\begin{definition}
Let $ (\mathcal{A}, \mu_{\mathcal{A}}, \alpha_{1}, \alpha_{2}) $ and $(\mathcal{B}, \mu_{\mathcal{B}}, \beta_{1}, \beta_{2}) $ 
be two partially hom-associative ternary algebras. Suppose that there are linear maps 
$\mathcal{L}_{\mu_{\mathcal{A}}}, \mathcal{M}_{\mu_{\mathcal{A}}}, \mathcal{R}_{\mu_{\mathcal{A}}}$ 
and $\mathcal{L}_{\mu_{\mathcal{B}}}, \mathcal{M}_{\mu_{\mathcal{B}}}, \mathcal{R}_{\mu_{\mathcal{B}}}$
 such that $(\mathcal{L}_{\mu_{\mathcal{A}}}, \mathcal{M}_{\mu_{\mathcal{A}}}, \mathcal{R}_{\mu_{\mathcal{A}}}, \beta_{1}, \beta_{2}, \mathcal{B})$ is a trimodule of $(\mathcal{A}, \mu_{\mathcal{A}}, \alpha_{1}, \alpha_{2}),$ 
and $(\mathcal{L}_{\mu_{\mathcal{B}}}, \mathcal{M}_{\mu_{\mathcal{B}}}, \mathcal{R}_{\mu_{\mathcal{B}}}, \alpha_{1}, \alpha_{2}, \mathcal{A})$ is a trimodule of $ (\mathcal{B}, \mu_{\mathcal{B}}, \beta_{1}, \beta_{2})$.Then 
$
(\mathcal{A}, \mathcal{B},\mathcal{L}_{\mu_{\mathcal{A}}}, \mathcal{M}_{\mu_{\mathcal{A}}}, \mathcal{R}_{\mu_{\mathcal{A}}}, \alpha_{1}, \alpha_{2}, \mathcal{L}_{\mu_{\mathcal{B}}}, \mathcal{M}_{\mu_{\mathcal{B}}}, \mathcal{R}_{\mu_{\mathcal{B}}}, \beta_{1}, \beta_{2})
$
 is a matched pair of partially hom-associative ternary algebras if Eqs.(\ref{pr7})-(\ref{pr26}) are satisfied, and if the following conditions hold:
 \beq \label{pr27}
 \mathcal{M}_{\mu_{A}}(\alpha_{1}(a), \alpha_{2}(z))(\mathcal{M}_{\mu_{A}}(\alpha_{1}(b), \alpha_{2}(y))(\mathcal{M}_{\mu_{A}}(\alpha_{1}(c), \alpha_{2}(x))\beta_{i}(v'))=\cr
  \mathcal{M}_{\mu_{A}}(\mu_{A}(\alpha_{1}(a), \alpha_{1}(b), \alpha_{1}(c)), \mu_{A}(\alpha_{2}(x), \alpha_{2}(y), \alpha_{2}(z))\beta_{i}(v'),
 \eeq   
 \beq \label{pr28}
 \mathcal{M}_{\mu_{B}}(\beta_{1}(a'), \beta_{2}(z'))(\mathcal{M}_{\mu_{B}}(\beta_{1}(b'), \beta_{2}(y'))(\mathcal{M}_{\mu_{B}}(\beta_{1}(c'), \beta_{2}(x'))\alpha_{i}(v))=\cr
  \mathcal{M}_{\mu_{B}}(\mu_{B}(\beta_{1}(a'), \beta_{1}(b'), \beta_{1}(c')), \mu_{B}(\beta_{2}(x'), \beta_{2}(y'), \beta_{2}(z')))\alpha_{i}(v),
 \eeq
 \beq\label{pr29}
\beta_{1}(\mathcal{L}_{\mu_{\mathcal{A}}}(a, b)v')= \mathcal{L}_{\mu_{\mathcal{A}}}(\alpha_{1}(a), \alpha_{2}(b))\beta_{1}(v'),\ \beta_{2}(\mathcal{L}_{\mu_{\mathcal{A}}}(a, b)v')= \mathcal{L}_{\mu_{\mathcal{A}}}(\alpha_{1}(a), \alpha_{2}(b))\beta_{2}(v')
\eeq
\beq\label{pr30}
\beta_{1}(\mathcal{M}_{\mu_{\mathcal{A}}}(a, b)v')= \mathcal{M}_{\mu_{\mathcal{A}}}(\alpha_{1}(a), \alpha_{2}(b))\beta_{1}(v'),\
\beta_{2}(\mathcal{M}_{\mu_{\mathcal{A}}}(a, b)v')= \mathcal{M}_{\mu_{\mathcal{A}}}(\alpha_{1}(a), \alpha_{2}(b))\beta_{2}(v')
\eeq
\beq\label{pr31}
\beta_{1}(\mathcal{R}_{\mu_{\mathcal{A}}}(a, b)v')= \mathcal{R}_{\mu_{\mathcal{A}}}(\alpha_{1}(a), \alpha_{2}(b))\beta_{1}(v'), 
\beta_{2}(\mathcal{R}_{\mu_{\mathcal{A}}}(a, b)v')= \mathcal{R}_{\mu_{\mathcal{A}}}(\alpha_{1}(a), \alpha_{2}(b))\beta_{2}(v')
\eeq
 \beq\label{pr32}
\alpha_{1}(\mathcal{L}_{\mu_{\mathcal{B}}}(a', b')v)= \mathcal{L}_{\mu_{\mathcal{B}}}(\beta_{1}(a'), \beta_{2}(b'))\alpha_{1}(v),\ \alpha_{2}(\mathcal{L}_{\mu_{\mathcal{B}}}(a', b')v)= \mathcal{L}_{\mu_{\mathcal{B}}}(\beta_{1}(a'), \beta_{2}(b'))\alpha_{2}(v)
\eeq
 \beq\label{pr33}
\alpha_{1}(\mathcal{M}_{\mu_{\mathcal{B}}}(a', b')v)= \mathcal{M}_{\mu_{\mathcal{B}}}(\beta_{1}(a'), \beta_{2}(b'))\alpha_{1}(v), \alpha_{2}(\mathcal{M}_{\mu_{\mathcal{B}}}(a', b')v)= \mathcal{M}_{\mu_{\mathcal{B}}}(\beta_{1}(a'), \beta_{2}(b'))\alpha_{2}(v)
\eeq
 \beq\label{pr34}
\alpha_{1}(\mathcal{R}_{\mu_{\mathcal{B}}}(a', b')v)= \mathcal{R}_{\mu_{\mathcal{B}}}(\beta_{1}(a'), \beta_{2}(b'))\alpha_{1}(v),\ \alpha_{2}(\mathcal{R}_{\mu_{\mathcal{B}}}(a', b')v)= \mathcal{R}_{\mu_{\mathcal{B}}}(\beta_{1}(a'), \beta_{2}(b'))\alpha_{2}(v)
\eeq
 for any $ x, y, z, a, b, c, v\in \mathcal{A}, x', y', z', a', b', c', v'\in \mathcal{B}$ and $i= 1, 2$.
\end{definition}

\section{Hom-associative ternary {infinitesimal} bialgebras}
We introduce here the notion of hom-associative ternary bialgebras,  characterized by a compatible condition between the ternary product and the ternary coproduct. We start with the following definitions.  
\begin{definition}\label{definition of PHB}
A partially hom-associative ternary {infinitesimal} bialgebra is a quadriple $(\mathcal{A}, \mu, \Delta, (\alpha_{i})_{i=1, 2})$ such that 
\begin{enumerate}
\item $(\mathcal{A}, \mu, (\alpha_{i})_{i=1, 2})$ is a partially hom-associative ternary algebra,
\item $(\mathcal{A}, \Delta, (\alpha_{i})_{i=1, 2})$ is a partially hom-coassociative ternary coalgebra,
\item $\Delta$ and $\mu$ satisfy the following condition:
\beq \textit{\label{cp1}}
\Delta\mu(x, y, z)&=&(L_{\mu}(\alpha_{1}(x), \alpha_{2}(y))\otimes\alpha_{1}\otimes\alpha_{2})\Delta(z) + (\alpha_{1}\otimes M_{\mu}(\alpha_{1}(x), \alpha_{2}(z))\otimes\alpha_{2})\Delta(y) +\cr &&(\alpha_{1}\otimes\alpha_{2}\otimes R_{\mu}(\alpha_{1}(y), \alpha_{2}(z))\Delta(x)
\eeq
where $ L_{\mu}(x, y)(z)= \mu(x, y, z); R_{\mu}(x, y)(z)= \mu(z, x, y); M_{\mu}(x, y)(z)= \mu(x, z, y)$ for $x, y,$ $ z \in \mathcal{A}.$
\end{enumerate} 
\end{definition}
\begin{example}\label{ep2}
Let $\mathcal{P}$ be a $2$-dimensionnal vector space with a basis $\lbrace e_{1}, e_{2}\rbrace$. Consider the following products defined by
\beqs
&& \mu: \mathcal{P}\otimes \mathcal{P}\otimes\mathcal{P}\rightarrow \mathcal{P}, \   \mu(e_{1}\otimes e_{1}\otimes e_{1})= e_{2}, \mu(e_{i}\otimes e_{j}\otimes e_{k})=0, i, j, k=1, 2,  (i, j, k)\neq (1, 1, 1); \cr
&& \Delta: \mathcal{P}\rightarrow \mathcal{P}\otimes \mathcal{P}\otimes\mathcal{P}, \  \Delta(e_{1})= e_{2}\otimes e_{2}\otimes e_{2}, \ \Delta(e_{2})= 0;\cr
&&\rho(e_{1})=  e_{1} +  e_{2},\ \rho(e_{2})=  e_{2},\ \rho= \alpha_{1}= \alpha_{2}.
\eeqs
By a straightforward  computation, one can easily show that the triple $(\mathcal{P}, \mu, \Delta, (\alpha_{i})_{i=1, 2})$ is a partially hom-associative ternary {infinitesimal} bialgebra.
\end{example} 

\begin{definition}\label{definition of HTB}
A totally hom-associative ternary {infinitesimal} bialgebra is a triple \\$(\mathcal{A}, \mu, \Delta, (\alpha_{i})_{i=1, 2})$ such that 
\begin{enumerate}
\item$(\mathcal{A}, \mu, (\alpha_{i})_{i=1, 2})$ is a totally hom-associative ternary algebra,
\item $(\mathcal{A}, \Delta, (\alpha_{i})_{i=1, 2})$ is a totally hom-coassociative ternary coalgebra,
\item$\Delta$ and $\mu$ satisfy the condition (\ref{cp1}).
\end{enumerate} 
\end{definition}
\begin{example}\label{et2}

Let $\mathcal{T}$ be a $2$-dimensionnal vector space with a basis $\lbrace e_{1}, e_{2}\rbrace$. Consider the following products defined by
\beqs
&&\mu(e_{1}\otimes e_{1}\otimes e_{1})= e_{1} \ \ \ \mu(e_{2}\otimes e_{2}\otimes e_{1})= 2e_{1} - e_{2}\cr
&&\mu(e_{1}\otimes e_{1}\otimes e_{2})= e_{1}-e_{2} \ \ \ \mu(e_{2}\otimes e_{2}\otimes e_{2})= 3e_{1} - 2e_{2}\cr
&&\mu(e_{1}\otimes e_{2}\otimes e_{2})= 2e_{1}-e_{2} \ \ \ \mu(e_{1}\otimes e_{2}\otimes e_{1})= e_{1} - e_{2}\cr
&&\mu(e_{2}\otimes e_{1}\otimes e_{1})= e_{1}-e_{2} \ \ \ \mu(e_{2}\otimes e_{1}\otimes e_{2})= 2e_{1} - e_{2},
\eeqs
\beqs
&& \Delta: \mathcal{T}\rightarrow \mathcal{T}\otimes \mathcal{T}\otimes\mathcal{T},\cr
\Delta(e_{1})&=& e_{1}\otimes e_{1}\otimes e_{1} + 2e_{1}\otimes e_{2}\otimes e_{2} + 2e_{2}\otimes e_{2}\otimes e_{1} + 3e_{2}\otimes e_{2}\otimes e_{2} + 2e_{2}\otimes e_{1}\otimes e_{2} +\cr && e_{1}\otimes e_{1}\otimes e_{2} + e_{2}\otimes e_{1}\otimes e_{1} + e_{1}\otimes e_{2}\otimes e_{1}\cr
\Delta(e_{2})&=& -e_{1}\otimes e_{1}\otimes e_{2} - e_{1}\otimes e_{2}\otimes e_{2} - e_{2}\otimes e_{1}\otimes e_{1} - e_{2}\otimes e_{2}\otimes e_{1} \cr
&& -2e_{2}\otimes e_{2}\otimes e_{2} -e_{1}\otimes e_{2}\otimes e_{1} -e_{2}\otimes e_{1}\otimes e_{2},
\eeqs
and the linear map
\beqs
\rho_{1}(e_{1})= e_{1},\
\rho_{1}(e_{2})= e_{1}-e_{2},\ \rho= \alpha_{1}= \alpha_{2}.
\eeqs
A direct computation shows that the triple $(\mathcal{T}, \mu, \Delta, (\alpha_{i})_{i=1, 2})$ is a totally hom-associative ternary {infinitesimal} bialgebra.
\end{example}
\begin{definition}
A weak totally hom-associative ternary {infinitesimal} bialgebra is a quadriple $(\mathcal{A}, \mu, \Delta, (\alpha_{i})_{i=1, 2})$ such that 
\begin{enumerate}
\item$(\mathcal{A}, \mu, (\alpha_{i})_{i=1, 2})$ is a weak totally hom-associative ternary algebra,
\item $(\mathcal{A}, \Delta, (\alpha_{i})_{i=1, 2})$ is a weak totally hom-coassociative ternary coalgebra,
\item $\Delta$ and $\mu$ satisfy the condition (\ref{cp1}).
\end{enumerate} 
\end{definition}

\begin{definition}
Two hom-associative ternary {infinitesimal} bialgebras $(\mathcal{A}_{1}, \mu_{1}, \Delta_{1}, (\alpha_{i})_{i=1, 2})$ and $(\mathcal{A}_{2}, \mu_{2}, \Delta_{2}, (\alpha'_{i})_{i=1, 2})$ are called equivalent if there exists a vector space isomorphism $f: \mathcal{A}_{1}\rightarrow \mathcal{A}_{2}$ such that
\begin{enumerate}
\item $f: (\mathcal{A}_{1}, \mu_{1}, (\alpha_{i})_{i=1, 2})\rightarrow (\mathcal{A}_{2}, \mu_{2}, (\alpha'_{i})_{i=1, 2})$ is a hom-associative ternary algebra isomorphism, that is, 
\beqs
f\mu_{1}(x, y, z)= \mu_{2}(f(x), f(y), f(z)) \mbox { for all } x, y, z \in \mathcal{A}_{1} \mbox{ and } f\circ \alpha_{i}=\alpha'_{i}\circ f \ \ \forall i= 1, 2.
\eeqs 
\item $f: (\mathcal{A}_{1}, \Delta_{1}, (\alpha_{i})_{i=1, 2})\rightarrow (\mathcal{A}_{2}, \Delta_{2}, (\alpha'_{i})_{i=1, 2})$ is a hom-coassociative ternary coalgebra isomorphism, that is, 
\beqs
\Delta_{2}(f(x))= (f\otimes f\otimes f)\Delta_{1}(x) \mbox { for every } x\in \mathcal{A}_{1} \mbox{ and } f\circ \alpha_{i}=\alpha'_{i}\circ f \ \ \forall i= 1, 2.
\eeqs 
\end{enumerate}
\end{definition} 

\begin{example}
Let $\mathcal{P}$ be a $2$-dimensionnal vector space with a basis $e_{1}, e_{2}$. Consider the following products $\mu_{1}, \mu_{2}: \mathcal{P}\otimes \mathcal{P}\otimes\mathcal{P}\rightarrow \mathcal{P},$\ $\Delta_{1}, \Delta_{2}: \mathcal{P}\rightarrow \mathcal{P}\otimes \mathcal{P}\otimes\mathcal{P},$ defined by
\beqs
&&\mu_{1}(e_{1}\otimes e_{1}\otimes e_{1})= e_{2}, \mu_{1}(e_{i}\otimes e_{j}\otimes e_{k})=0, i, j, k=1, 2,  (i, j, k)\neq (1, 1, 1), \cr
&&\Delta_{1}(e_{1})= e_{2}\otimes e_{2}\otimes e_{2}, \ \Delta_{1}(e_{2})= 0,\cr
&&\rho(e_{1})=  e_{1} +  e_{2},\ \rho(e_{2})=  e_{2},\ \rho= \alpha_{1}= \alpha_{2},
\eeqs 
and
\beqs
&&\mu_{2}(e_{2}\otimes e_{2}\otimes e_{2})= e_{1}, \mu_{2}(e_{i}\otimes e_{j}\otimes e_{k})=0, i, j, k=1, 2,  (i, j, k)\neq (2, 2, 2), \cr
&&\Delta_{2}(e_{2})= e_{1}\otimes e_{1}\otimes e_{1}, \ \Delta_{2}(e_{1})= 0,\cr
&&\rho'(e_{2})=  e_{1} +  e_{2},\ \rho'(e_{1})=  e_{1},\ \rho'= \alpha'_{1}= \alpha'_{2}.
\eeqs 
By a direct computation, we can prove that  the triples $(\mathcal{P}, \mu_{1}, \Delta_{1}, (\alpha_{i})_{i=1, 2})$ and $(\mathcal{P}, \mu_{2}, \Delta_{2}, (\alpha'_{i})_{i=1, 2})$ are partially hom-associative ternary {infinitesimal} bialgebras. Let us consider the linear map $f: \mathcal{P}\rightarrow \mathcal{P}, f(e_{1})= e_{2}, f(e_{2})= e_{1}$. Then, $(\mathcal{P}, \mu_{1}, \Delta_{1}, (\alpha_{i})_{i=1, 2})$ and $(\mathcal{P}, \mu_{2}, \Delta_{2}, (\alpha'_{i})_{i=1, 2})$ are equivalent partially hom-associative ternary {infinitesimal} bialgebras by the isomorphism $f:(\mathcal{P}, \mu_{1}, \Delta_{1}, (\alpha_{i})_{i=1, 2})\rightarrow (\mathcal{P}, \mu_{2}, \Delta_{2}, (\alpha'_{i})_{i=1, 2})$.
\end{example}
{
Let $\A$ be a hom-associative ternary algebra, {and} $\sigma: \A\otimes \A\rightarrow \A\otimes\A$ {an} exchange operator defined as
\beqs
\sigma(a\otimes b)= b\otimes a, \mbox{ for all } a, b\in \A.
\eeqs
The condition (\ref{cp1}) can be rewritten as
\beq \label{IHATB}
\Delta\mu=(\mu\otimes\alpha_{1}\otimes\alpha_{2})(\alpha_{1}\otimes\alpha_{2}\otimes\Delta) + (\alpha_{1}\otimes\mu\otimes\alpha_{2})(\sigma\otimes \id\otimes\sigma)(\alpha_{1}\otimes\Delta\otimes\alpha_{2})\cr + (\alpha_{1}\otimes\alpha_{2}\otimes\mu)(\Delta\otimes\alpha_{1}\otimes\alpha_{2}).
\eeq 
An associative ternary infinitesimal bialgebra  is a hom-associative ternary infinitesimal bialgebra with $\alpha_{1}=\alpha_{2}=\id$.

In terms of elements and writing $\mu(a, b, c)$ as $abc$, the condition (\ref{IHATB}) can be rewritten as
\beq
\Delta(abc)=\alpha_{1}(a)\alpha_{2}(b)c_{1}\otimes\alpha(c_{2})\otimes\alpha(c_{3}) + \alpha_{1}(b_{1})\otimes\alpha_{1}(a)b_{2}\alpha_{2}(c)\otimes\alpha_{2}(b_{3})\cr + \alpha_{1}(a_{1})\otimes\alpha_{2}(a_{2})\otimes a_{3}\alpha_{1}(b)\alpha_{2}(c), \mbox{ for all } a, b, c\in \A.
\eeq
To be more explicit, an associative ternary infinitesimal bialgebra is a triple $(\A, \mu, \Delta)$ consisting of an associative ternary algebra $(\A, \mu)$ and a coassociative ternary coalgebra $(\A, \Delta)$ such that the condition
\beq \label{IATB}
\Delta\mu=(\mu\otimes\id\otimes\id)(\id\otimes\id\otimes\Delta) + (\id\otimes\mu\otimes\id)(\sigma\otimes \id\otimes\sigma)(\id\otimes\Delta\otimes\id)\cr + (\id\otimes\id\otimes\mu)(\Delta\otimes\id\otimes\id),
\eeq
{or, equivalently,}
\beq
\Delta(abc)=abc_{1}\otimes c_{2}\otimes c_{3} + b_{1}\otimes ab_{2}c\otimes b_{3} + a_{1}\otimes a_{2}\otimes a_{3}bc, \mbox{ for all } a, b, c\in \A
\eeq
{holds}.
\begin{example}
A hom-associative ternary algebra $(\A, \mu, \alpha_{1}, \alpha_{2})$ becomes a hom-associative 
ternary infinitesimal bialgebra when equipped with the trivial comultiplication $\Delta = 0$. 
Likewise, a hom-coassociative ternary coalgebra $(C, \Delta, \alpha)$ becomes a hom-associative 
ternary infinitesimal bialgebra when equipped with the trivial multiplication $\mu = 0$. 
\end{example}

 The dual maps $L^{*}_{\mu}, R^{*}_{\mu}, M^{*}_{\mu}$  of the linear maps $L_{\mu}, R_{\mu}, M_{\mu}$  (see Example \ref{extha}) are defined, respectively, as:\\ 
$\displaystyle L^{*}_{\mu}, R^{*}_{\mu}, M^{*}_{\mu}: \A\otimes \A \rightarrow \mathfrak{gl}(\A^{*})$
 such that:
\beq\label{dual1}
 L^*_{\mu}: \A\otimes \A & \longrightarrow & \mathfrak{gl}(\A^*)  \cr
  (x, y)  & \longmapsto & L^*_\mu(x, y):
      \begin{array}{llll}
 \A^* &\longrightarrow & \A^* \\ 
  u^* & \longmapsto & L^*_\mu(x, y)(u^*): 
      \begin{array}{llll}
\A  \longrightarrow   \mathcal{K} \cr
v  \longmapsto  &\left< L^{*}_\mu(x, y)(u^{*}), v\right>\cr
 &:= \left<  u^{*}, L_\mu(x, y) v\right>, 
\end{array} 
     \end{array}
 \eeq
\beq\label{dual2}
 R^*_{\mu}: \A\otimes \A & \longrightarrow & \mathfrak{gl}(\A^*)  \cr
  (x, y)  & \longmapsto & R^*_\mu(x, y):
      \begin{array}{llll}
 \A^* &\longrightarrow & \A^* \\ 
  u^* & \longmapsto & R^*_\mu(x, y)(u^*): 
      \begin{array}{llll}
\A  \longrightarrow  \mathcal{K} \cr
v  \longmapsto & \left< R^{*}_\mu(x, y)(u^{*}), v\right>\cr
 &:= \left<  u^{*}, R_\mu(x, y) v\right>, 
\end{array} 
     \end{array}
     \eeq
 
  \beq\label{dual3}
 M^*_{\mu}: \A\otimes \A & \longrightarrow & \mathfrak{gl}(\A^*)  \cr
  (x, y)  & \longmapsto & M^*_\mu(x, y):
      \begin{array}{llll}
 \A^* &\longrightarrow & \A^* \\ 
  u^* & \longmapsto & M^*_\mu(x, y)(u^*): 
      \begin{array}{llll}
\A  \longrightarrow  \mathcal{K} \cr
v\longmapsto &\left< M^{*}_\mu(x, y)(u^{*}), v\right>\cr
 &:= \left<  u^{*}, M_\mu(x, y) v\right>, 
\end{array} 
     \end{array}
 \eeq
for all $x, y, z, v \in \A, u^{*} \in \A^{*},$ where $\A^{*}$ is the dual space of $\A.$ 
}
\begin{theorem}\label{pb1}
Let $(\mathcal{P}, \mu, \Delta, (\alpha_{i})_{i=1, 2})$ be a partially hom-associative ternary infinitesimal bialgebra. Then, 
$(\mathcal{P}^{\ast}, \Delta^{\ast}, \mu^{\ast}, (\alpha^{\ast}_{i})_{i=1, 2})$ is a partially  hom-associative ternary infinitesimal bialgebra, called the dual partially hom-associative ternary infinitesimal bialgebra of $(\mathcal{P}, \mu, \Delta, (\alpha_{i})_{i=1, 2})$.
\end{theorem}
\textbf{Proof:}
$(\mathcal{P}, \mu, \Delta, (\alpha_{i})_{i=1, 2})$ is a partially hom-associative ternary infinitesimal bialgebra. Then, by the Theorem (\ref{teo-co-asso}) and the Theorem (\ref{teo-co-asso1}), $(\mathcal{P}^{\ast}, \Delta^{\ast}, (\alpha^{\ast}_{i})_{i=1, 2})$ is a partially hom-associative ternary algebra with the multiplication (\ref{eqco}), and $(\mathcal{P}^{\ast}, \mu^{\ast}, (\alpha^{\ast}_{i})_{i=1, 2})$ is a partially hom-coassociative ternary coalgebra with the multiplication (\ref{qdual}). Now, we prove that $\mu^{\ast}: \mathcal{P}^{\ast}\rightarrow \mathcal{P}^{\ast}\otimes \mathcal{P}^{\ast}\otimes \mathcal{P}^{\ast}$ satisfies the identity (\ref{cp1}), that is, the following identity holds for every $\xi, \eta, \gamma \in \mathcal{P}^{\ast} $,
 \beq \textit{\label{cp2}}
\mu^{\ast}(\Delta^{\ast}(\xi, \eta, \gamma))&=&(L_{\Delta^{\ast}}(\alpha^{\ast}_{1}(\xi), \alpha^{\ast}_{2}(\eta))\otimes\alpha^{\ast}_{1}\otimes\alpha^{\ast}_{2})\mu^{\ast}(\gamma) + (\alpha^{\ast}_{1}\otimes M_{\Delta^{\ast}}(\alpha^{\ast}_{1}(\xi), \alpha^{\ast}_{2}(\gamma))\otimes\alpha^{\ast}_{2})\mu^{\ast}(\eta)\cr
&& + (\alpha^{\ast}_{1}\otimes\alpha^{\ast}_{2}\otimes R_{\Delta^{\ast}}(\alpha^{\ast}_{1}(\eta), \alpha^{\ast}_{2}(\gamma))\mu^{\ast}(\xi).
\eeq
For every $x, y, z\in \mathcal{P}$ and $\xi, \eta, \gamma \in \mathcal{P}^{\ast},$ by identities (\ref{eqco}) and (\ref{qdual}),
\beqs
\langle \mu^{\ast}\Delta^{\ast}(\xi, \eta, \gamma), x\otimes y\otimes z \rangle 
&=&\langle \xi\otimes \eta\otimes \gamma, \Delta\mu(x\otimes y\otimes z)  \rangle \cr
&=&\langle \xi\otimes \eta\otimes \gamma, (L_{\mu}(\alpha_{1}(x), \alpha_{2}(y))\otimes\alpha_{1}\otimes\alpha_{2})\Delta(z) +\cr
&& (\alpha_{1}\otimes M_{\mu}(\alpha_{1}(x), \alpha_{2}(z))\otimes\alpha_{2})\Delta(y) +\cr
&& (\alpha_{1}\otimes\alpha_{2}\otimes R_{\mu}(\alpha_{1}(y), \alpha_{1}(z)))\Delta(x)\rangle\cr
&=&\langle \xi\otimes \eta\otimes \gamma, (L_{\mu}(\alpha_{1}(x), \alpha_{2}(y))\otimes\alpha_{1}\otimes\alpha_{2})\Delta(z)\rangle +\cr 
&&\langle\xi\otimes \eta\otimes \gamma, (\alpha_{1}\otimes M_{\mu}(\alpha_{1}(x), \alpha_{2}(z))\otimes\alpha_{2})\Delta(y)\rangle +\cr
&&\langle \xi\otimes \eta\otimes \gamma, (\alpha_{1}\otimes\alpha_{2}\otimes R_{\mu}(\alpha_{1}(y), \alpha_{2}(z)))\Delta(x)\rangle \cr
&=&\langle \xi\otimes \eta\otimes \gamma, ([(L_{\mu}\circ(\alpha_{1}, \alpha_{2}))(x,y)]\otimes\alpha_{1}\otimes\alpha_{2})\Delta(z)\rangle +\cr
&&\langle\xi\otimes \eta\otimes \gamma, (\alpha_{1}\otimes [(M_{\mu}\circ(\alpha_{1}, \alpha_{2}))(x, z)]\otimes\alpha_{2})\Delta(y)\rangle +\cr
&&\langle \xi\otimes \eta\otimes \gamma, (\alpha_{1}\otimes\alpha_{2}\otimes [(R_{\mu}\circ(\alpha_{1}, \alpha_{2}))(y, z)])\Delta(x)\rangle \cr
&=&\langle \Delta^{\ast}(\alpha^{\ast}_{1}\otimes\alpha^{\ast}_{2}\otimes[((\alpha^{\ast}_{1}, \alpha^{\ast}_{2})\circ L^{\ast}_{\mu})(x,y)])(\xi\otimes \eta\otimes \gamma), z\rangle +\cr
&& \langle \Delta^{\ast}(\alpha^{\ast}_{1}\otimes [((\alpha^{\ast}_{1}, \alpha^{\ast}_{2})\circ M^{\ast}_{\mu})(x,z)]\otimes\alpha^{\ast}_{2})(\xi\otimes \eta\otimes \gamma), y\rangle +\cr
&&  \langle\Delta^{\ast}([((\alpha^{\ast}_{1}, \alpha^{\ast}_{2})\circ R^{\ast}_{\mu})(y, z)]\otimes\alpha^{\ast}_{1}\otimes\alpha^{\ast}_{2})(\xi\otimes \eta\otimes \gamma), x\rangle\cr
&=&\langle \Delta^{\ast}(\alpha^{\ast}_{1}(\xi)\otimes \alpha^{\ast}_{2}(\eta)\otimes[((\alpha^{\ast}_{1}, \alpha^{\ast}_{2})\circ L^{\ast}_{\mu})(x,y)](\gamma)), z\rangle + \cr
&&\langle \Delta^{\ast}(\alpha^{\ast}_{1}(\xi)\otimes [((\alpha^{\ast}_{1}, \alpha^{\ast}_{2})\circ M^{\ast}_{\mu})(x,z)](\eta)\otimes \alpha_{2}(\gamma)), y\rangle +\cr
 &&\langle\Delta^{\ast}([((\alpha^{\ast}_{1}, \alpha^{\ast}_{2})\circ R^{\ast}_{\mu})(y, z)](\xi)\otimes\alpha^{\ast}_{1}(\eta)\otimes \alpha^{\ast}_{2}(\gamma)), x\rangle\cr
&=&\langle L_{\Delta^{\ast}}(\alpha^{\ast}_{1}(\xi), \alpha^{\ast}_{2}(\eta))([((\alpha^{\ast}_{1}, \alpha^{\ast}_{2})\circ L^{\ast}_{\mu})(x,y)](\gamma)), z\rangle +\cr
&&\langle M_{\Delta^{\ast}}(\alpha^{\ast}_{1}(\xi), \alpha^{\ast}_{2}(\gamma))([((\alpha^{\ast}_{1}, \alpha^{\ast}_{2})\circ M^{\ast}_{\mu})(x,z)](\eta)), y\rangle +\cr
&& \langle R_{\Delta^{\ast}}( \alpha^{\ast}_{1}(\eta), \alpha^{\ast}_{2}(\gamma))([((\alpha^{\ast}_{1}, \alpha^{\ast}_{2})\circ R^{\ast}_{\mu})(y, z)](\xi)), x\rangle\cr
&=&\langle [((\alpha^{\ast}_{1}, \alpha^{\ast}_{2})\circ L^{\ast}_{\mu})(x,y)](\gamma), L^{\ast}_{\Delta^{\ast}}(\alpha^{\ast}_{1}(\xi), \alpha^{\ast}_{2}(\eta))(z)\rangle +\cr
&&\langle [((\alpha^{\ast}_{1}, \alpha^{\ast}_{2})\circ M^{\ast}_{\mu})(x,z)](\eta), M^{\ast}_{\Delta^{\ast}}(\alpha^{\ast}_{1}(\xi), \alpha^{\ast}_{2}(\gamma))(y)\rangle +\cr
&& \langle [((\alpha^{\ast}_{1}, \alpha^{\ast}_{2})\circ R^{\ast}_{\mu})(y, z)](\xi), R^{\ast}_{\Delta^{\ast}}(\alpha^{\ast}_{1}(\eta), \alpha^{\ast}_{2}(\gamma))(x)\rangle\cr
&=&\langle \gamma, \mu(\alpha_{1}(x), \alpha_{2}(y), L^{\ast}_{\Delta^{\ast}}(\alpha^{\ast}_{1}(\xi), \alpha^{\ast}_{2}(\eta))(z))\rangle +\cr
&& \langle \eta, \mu(\alpha_{1}(x), M^{\ast}_{\Delta^{\ast}}(\alpha^{\ast}_{1}(\xi), \alpha^{\ast}_{2}(\gamma))(y), \alpha_{2}(z)\rangle +\cr
&& \langle \xi, \mu(R^{\ast}_{\Delta^{\ast}}( \alpha^{\ast}_{1}(\eta), \alpha^{\ast}_{2}(\gamma))(x), \alpha_{1}(y), \alpha_{2}(z))\rangle \cr
&=&\langle \mu^{\ast}(\gamma), \alpha_{1}(x)\otimes \alpha_{2}(y)\otimes L^{\ast}_{\Delta^{\ast}}(\alpha^{\ast}_{1}(\xi), \alpha^{\ast}_{2}(\eta))(z))\rangle +\cr
&& \langle \mu^{\ast}(\eta), \alpha_{1}(x)\otimes M^{\ast}_{\Delta^{\ast}}(\alpha^{\ast}_{1}(\xi), \alpha^{\ast}_{2}(\gamma))(y)\otimes \alpha_{2}(z)\rangle +\cr
&& \langle \mu^{\ast}(\xi), R^{\ast}_{\Delta^{\ast}}( \alpha^{\ast}_{1}(\eta), \alpha^{\ast}_{2}(\gamma))(x)\otimes \alpha_{1}(y)\otimes \alpha_{2}(z)\rangle\cr
&=&\langle \mu^{\ast}(\gamma), (\alpha_{1}\otimes \alpha_{2}\otimes L^{\ast}_{\Delta^{\ast}}(\alpha^{\ast}_{1}(\xi), \alpha^{\ast}_{2}(\eta)))(x\otimes y\otimes z)\rangle +\cr
&& \langle \mu^{\ast}(\eta), (\alpha_{1}\otimes M^{\ast}_{\Delta^{\ast}}(\alpha^{\ast}_{1}(\xi), \alpha^{\ast}_{2}(\gamma))\otimes \alpha_{2})(x\otimes y\otimes z)\rangle +\cr
&& \langle \mu^{\ast}(\xi), (R^{\ast}_{\Delta^{\ast}}(\alpha^{\ast}_{1}(\eta), \alpha^{\ast}_{2}(\gamma))\otimes \alpha_{1}\otimes \alpha_{2})(x\otimes y\otimes z)\rangle\cr
&=&\langle (\alpha^{\ast}_{1}\otimes \alpha^{\ast}_{2}\otimes R_{\Delta^{\ast}}(\alpha^{\ast}_{1}(\eta), \alpha^{\ast}_{2}(\gamma)))\mu^{\ast}(\xi), x\otimes y\otimes z\rangle + \cr
&&\langle (\alpha^{\ast}_{1}\otimes M_{\Delta^{\ast}}(\alpha^{\ast}_{1}(\xi), \alpha^{\ast}_{2}(\gamma))\otimes \alpha^{\ast}_{2})\mu^{\ast}(\eta), x\otimes y\otimes z\rangle +\cr
&& \langle(L_{\Delta^{\ast}}(\alpha^{\ast}_{1}(\xi), \alpha^{\ast}_{2}(\eta))\otimes \alpha^{\ast}_{1}\otimes \alpha^{\ast}_{2})\mu^{\ast}(\gamma), x\otimes y\otimes z\rangle.
\eeqs
Then, we obtain
\beqs
 \mu^{\ast}\Delta^{\ast}(\xi, \eta, \gamma)&=& (\alpha^{\ast}_{1}\otimes \alpha^{\ast}_{2}\otimes R_{\Delta^{\ast}}(\alpha^{\ast}_{1}(\eta), \alpha^{\ast}_{2}(\gamma)))\mu^{\ast}(\xi)+ (\alpha^{\ast}_{1}\otimes M_{\Delta^{\ast}}(\alpha^{\ast}_{1}(\xi), \alpha^{\ast}_{2}(\gamma))\otimes \alpha^{\ast}_{2})\mu^{\ast}(\eta)+\cr
 &&(L_{\Delta^{\ast}}(\alpha^{\ast}_{1}(\xi), \alpha^{\ast}_{2}(\eta))\otimes \alpha^{\ast}_{1}\otimes \alpha^{\ast}_{2})\mu^{\ast}(\gamma).
\eeqs
Hence, the identity (\ref{cp2}) holds.

$ \hfill \square $
\begin{example}
From Theorem (\ref{pb1}), the dual partially hom-associative ternary infinitesimal bialgebra ($\mathcal{P}^{\ast}, \Delta^{\ast}, \mu^{\ast}, (\alpha^{\ast}_{i})_{i=1, 2}$) of ($\mathcal{P}, \mu, \Delta, (\alpha_{i})_{i=1, 2}$) in Example (\ref{ep2}), has its multiplications\\ $\Delta^{\ast}: \mathcal{P}^{\ast}\otimes \mathcal{P}^{\ast}\otimes\mathcal{P}^{\ast}\rightarrow \mathcal{P}^{\ast}$ and 
$ \mu^{\ast}: \mathcal{P}^{\ast}\rightarrow \mathcal{P}^{\ast}\otimes \mathcal{P}^{\ast}\otimes\mathcal{P}^{\ast},$ defined by:
\beqs
&&\Delta^{\ast}(e^{\ast}_{2}\otimes e^{\ast}_{2}\otimes e^{\ast}_{2})= e^{\ast}_{1}, \Delta^{\ast}(e^{\ast}_{i}\otimes e^{\ast}_{j}\otimes e^{\ast}_{k})=0, i, j, k=1, 2, (i, j, k)\neq (2, 2, 2); \cr
&&\mu^{\ast}(e^{\ast}_{2})= e^{\ast}_{1}\otimes e^{\ast}_{1}\otimes e^{\ast}_{1}, \mu^{\ast}(e^{\ast}_{1})=0;\cr
&&\rho^{\ast}(e^{\ast}_{1})=  e^{\ast}_{1},\ \rho^{\ast}(e^{\ast}_{2})= e^{\ast}_{1} + e^{\ast}_{2},\ \rho^{\ast}= \alpha^{\ast}_{1}= \alpha^{\ast}_{2}.
\eeqs 
\end{example}

Similarly, we have the following results.
\begin{theorem}\label{tb1}
Let $(\mathcal{A}, \mu, \Delta, (\alpha_{i})_{i=1, 2})$ be a totally hom-associative ternary infinitesimal bialgebra. Then, 
$(\mathcal{A}^{\ast}, \Delta^{\ast}, \mu^{\ast}, (\alpha^{\ast}_{i})_{i=1, 2})$ is a totally hom-associative ternary infinitesimal bialgebra,  called the dual totally hom-associative ternary infinitesimal bialgebra of $(\mathcal{A}, \mu, \Delta, (\alpha_{i})_{i=1, 2})$.
\end{theorem}
\begin{example}
From Theorem (\ref{tb1}), the dual totally hom-associative ternary infinitesimal bialgebra ($\mathcal{T}^{\ast}, \Delta^{\ast}, \mu^{\ast}, \rho^{\ast}$) of ($\mathcal{T}, \mu, \Delta, \rho$) in Example (\ref{et2}), has its multiplications $\Delta^{\ast}: \mathcal{T}^{\ast}\otimes \mathcal{T}^{\ast}\otimes\mathcal{T}^{\ast}\rightarrow \mathcal{T}^{\ast}$ and 

 $ \mu^{\ast}: \mathcal{T}^{\ast}\rightarrow \mathcal{T}^{\ast}\otimes \mathcal{T}^{\ast}\otimes\mathcal{T}^{\ast},$ defined such that $\Delta^{\ast}= \mu$ and $\mu^{\ast}= \Delta$.
\end{example} 
\begin{theorem}
Let $(\mathcal{A}, \mu, \Delta, (\alpha_{i})_{i=1, 2})$ be a weak totally hom-associative ternary infinitesimal bialgebra. Then, 
$(\mathcal{A}^{\ast}, \Delta^{\ast}, \mu^{\ast}, (\alpha^{\ast}_{i})_{i=1, 2})$ is a weak totally hom-associative ternary infinitesimal bialgebra, called the dual weak totally hom-associative ternary infinitesimal bialgebra of $(\mathcal{A}, \mu, \Delta)$.
\end{theorem}
Next, let us characterize  the hom-associative ternary infinitesimal bialgebras in terms of {algebra structure} constants.
For that, consider a partially hom-associative ternary infinitesimal bialgebra  $(\mathcal{A}, \mu, \Delta, (\alpha_{i})_{i=1, 2})$ with the multiplications in the basis $e_{1},..., e_{n}$ defined as follows:
\beq \label{cp3}
\mu(e_{i}, e_{j}, e_{k})=\sum^{n}_{l=1} c^{l}_{ijk}e_{l},\ \Delta(e_{l})=\sum_{1\leq r, s, t\leq n}a^{rst}_{l}e_{r}\otimes e_{s}\otimes e_{t},\ \alpha_{q}(e_{j})=\sum^{n}_{k=1}y^{qj}_{k}e_{k},\ q=1, 2;
\eeq
where $c^{l}_{ijk}, a^{rst}_{l}\in \mathcal{K}, 1\leq i,j, k, l\leq n.$ Then,
\beq \label{cp5}
\Delta\mu(e_{i}, e_{j}, e_{k})=\Delta(\sum^{n}_{l=1} c^{l}_{ijk}e_{l})= \sum^{n}_{l=1} c^{l}_{ijk}\Delta(e_{l})= \sum^{n}_{l=1}\sum_{1\leq r, s, t\leq n}c^{l}_{ijk} a^{rst}_{l}e_{r}\otimes e_{s}\otimes e_{t} 
\eeq
From the identity (\ref{cp1}), we derive
\beq \label{cp4}
\Delta\mu(e_{i}, e_{j}, e_{k})&=&(L_{\mu}(\alpha_{1}(e_{i}), \alpha_{2}(e_{j}))\otimes\alpha_{1}\otimes\alpha_{2})\Delta(e_{k}) +\cr
&&(\alpha_{1}\otimes M_{\mu}(\alpha_{1}(e_{i}), \alpha_{2}(e_{k}))\otimes\alpha_{2})\Delta(e_{j}) +\cr
&& (\alpha_{1}\otimes\alpha_{2}\otimes R_{\mu}(\alpha_{1}(e_{j}), \alpha_{2}(e_{k})))\Delta(e_{i})\cr
&=&(L_{\mu}(\alpha_{1}(e_{i}), \alpha_{2}(e_{j}))\otimes\alpha_{1}\otimes\alpha_{2})\left( \sum_{1\leq r, s, t\leq n}a^{rst}_{k}e_{r}\otimes e_{s}\otimes e_{t}\right)  +\cr
&& (\alpha_{1}\otimes M_{\mu}(\alpha_{1}(e_{i}), \alpha_{2}(e_{k}))\otimes\alpha_{2})\left( \sum_{1\leq r, s, t\leq n}a^{rst}_{j}e_{r}\otimes e_{s}\otimes e_{t}\right)  +\cr
&&(\alpha_{1}\otimes\alpha_{2}\otimes R_{\mu}(\alpha_{1}(e_{j}), \alpha_{2}(e_{k})))\left( \sum_{1\leq r, s, t\leq n}a^{rst}_{i}e_{r}\otimes e_{s}\otimes e_{t}\right) \cr
&=&(L_{\mu}\left(\sum^{n}_{p=1}y^{1i}_{p}e_{p}, \sum^{n}_{q=1}y^{2j}_{q}e_{q}\right) \otimes\alpha_{1}\otimes\alpha_{2})\left( \sum_{1\leq r, s, t\leq n}a^{rst}_{k}e_{r}\otimes e_{s}\otimes e_{t}\right)+\cr
&& (\alpha_{1}\otimes M_{\mu}\left(\sum^{n}_{p=1}y^{1i}_{p}e_{p}, \sum^{n}_{q=1}y^{2k}_{q}e_{q}\right)\otimes\alpha_{2})\left( \sum_{1\leq r, s, t\leq n}a^{rst}_{j}e_{r}\otimes e_{s}\otimes e_{t}\right)  +\cr
&&(\alpha_{1}\otimes\alpha_{2}\otimes R_{\mu}\left(\sum^{n}_{p=1}y^{1j}_{p}e_{p}, \sum^{n}_{q=1}y^{2k}_{q}e_{q}\right))\left(\sum_{1\leq r, s, t\leq n}a^{rst}_{i}e_{r}\otimes e_{s}\otimes e_{t}\right) \cr
&=&\sum_{1\leq r, s, t, p, q\leq n}a^{rst}_{k}y^{1i}_{p}y^{2j}_{q}\mu(e_{p}, e_{q}, e_{r})\otimes \alpha_{1}(e_{s})\otimes \alpha_{2}(e_{t})+\cr
 && \sum_{1\leq r, s, t, p, q\leq n}a^{rst}_{j}y^{1i}_{p}y^{2k}_{q}\alpha_{1}(e_{r})\otimes \mu(e_{p}, e_{s}, e_{q})\otimes \alpha_{2}(e_{t}) +\cr
&& \sum_{1\leq r, s, t, p, q\leq n}a^{rst}_{i}y^{1j}_{p}y^{2k}_{q}\alpha_{1}(e_{r})\otimes \alpha_{2}(e_{s})\otimes \mu(e_{t},e_{p}, e_{q}) \cr
&=&\sum_{1\leq r, s, t, p, q\leq n}a^{rst}_{k}y^{1i}_{p}y^{2j}_{q}\left(\sum^{n}_{l=1} c^{l}_{pqr}e_{l}\right) \otimes\left(\sum^{n}_{u=1}y^{1s}_{u}e_{u}\right) \otimes \left(\sum^{n}_{v=1}y^{2t}_{v}e_{v}\right)+\cr
&&
\sum_{1\leq r, s, t, p, q\leq    n}a^{rst}_{j}y^{1i}_{p}y^{2k}_{q}\left(\sum^{n}_{u=1}y^{1r}_{u}e_{u}\right)\otimes \left(\sum^{n}_{l=1} c^{l}_{psq}e_{l}\right)\otimes\left(\sum^{n}_{v=1}y^{2t}_{v}e_{v}\right) +\cr
&& \sum_{1\leq r, s, t, p, q\leq n}a^{rst}_{i}y^{1j}_{p}y^{2k}_{q}\left(\sum^{n}_{u=1}y^{1r}_{u}e_{u}\right)\otimes \left(\sum^{n}_{v=1}y^{2s}_{v}e_{v}\right)\otimes \left(\sum^{n}_{l=1} c^{l}_{tpq}e_{l}\right)\cr
&=&\sum_{1\leq r, s, t, p, q, l, u, v\leq n}a^{rst}_{k}c^{l}_{pqr}y^{1i}_{p}y^{2j}_{q}y^{1s}_{u}y^{2t}_{v} e_{l} \otimes e_{u} \otimes e_{v}+\cr
&&
\sum_{1\leq r, s, t, p, q, l, u, v\leq n}a^{rst}_{j}c^{l}_{psq}y^{1i}_{p}y^{2k}_{q}y^{1r}_{u} y^{2t}_{v} e_{u}\otimes e_{l}\otimes e_{v} +\cr
&& \sum_{1\leq r, s, t, p, q, l, u, v\leq n}a^{rst}_{i} c^{l}_{tpq}y^{1j}_{p}y^{2k}_{q} y^{1r}_{u}y^{2s}_{v} e_{u}\otimes e_{v}\otimes e_{l}\cr
&=&\sum_{1\leq r, s, t, p, q, l, u, v\leq n}[a^{rst}_{k}c^{l}_{pqr}y^{1i}_{p}y^{2j}_{q}y^{1s}_{u}y^{2t}_{v}+ a^{rst}_{j}c^{u}_{psq}y^{1i}_{p}y^{2k}_{q}y^{1r}_{l} y^{2t}_{v} +\cr
&& a^{rst}_{i} c^{v}_{tpq}y^{1j}_{p}y^{2k}_{q} y^{1r}_{l}y^{2s}_{u}]e_{l} \otimes e_{u} \otimes e_{v}.
\eeq
Comparing the identities (\ref{cp5}) and (\ref{cp4}), we arrive at:
\beq
\sum_{l=1}^{n}[a^{rst}_{k}c^{l}_{pqr}y^{1i}_{p}y^{2j}_{q}y^{1s}_{u}y^{2t}_{v}+ a^{rst}_{j}c^{u}_{psq}y^{1i}_{p}y^{2k}_{q}y^{1r}_{l} y^{2t}_{v} + a^{rst}_{i} c^{v}_{tpq}y^{1j}_{p}y^{2k}_{q} y^{1r}_{l}y^{2s}_{u} -a^{rst}_{l}c^{l}_{ijk}]=0.
\eeq
Conversely, if a hom-associative ternary algebra ($\mathcal{A}, \mu, (\alpha_{i})_{i=1, 2}$) and a hom-coassociative ternary coalgebra ($\mathcal{A}, \Delta, (\alpha_{i})_{i=1, 2}$) defined by (\ref{cp3}), satisfy the identity (\ref{cp4}), then $\mu$ and $\Delta$ obey the identity (\ref{cp1}). Therefore, we conclude by the following fundamental result:

\begin{theorem}
Let $\mathcal{A}$ be a vector space with a basis $e_{1},..., e_{n}$ and $(\mathcal{A}, \mu, (\alpha_{i})_{i=1, 2})$ and $(\mathcal{A}, \Delta, (\alpha_{i})_{i=1, 2})$ be hom-associative ternary algebra and hom-coassociative ternary coalgebra defined by (\ref{cp3}). Then, $(\mathcal{A}, \mu, \Delta, (\alpha_{i})_{i=1, 2})$ is an hom-associative ternary infinitesimal bialgebra if and only if $c^{l}_{ijk}$ and $a^{rst}_{l},$ $1\leq i, j, k, l, r, s, t\leq n,$ satisfy {the} identity (\ref{cp4}). 
\end{theorem}
{
\begin{definition}\label{definition of inf. ternary hom-bialgebra}
A hom-associative ternary infinitesimal bialgebra is a triple \\$(\A, \mu, \Delta, (\alpha_{i})_{i=1, 2}))$ consisting of a hom-associative ternary algebra $(\A, \mu, (\alpha_{i})_{i=1, 2}))$ and a hom-coassociative ternary coalgebra $(\A, \Delta, (\alpha_{i})_{i=1, 2}))$ such that the condition (\ref{IATB}) holds.
\end{definition}
\begin{proposition}\label{relevent proposition}
Let $(\A, \mu, \Delta, \alpha_{1}, \alpha_{2})$ be a hom-associative ternary infinitesimal
bialgebra. Then so are $(\A, -\mu, \Delta, \alpha_{1}, \alpha_{2}), (\A, \mu, -\Delta, \alpha_{1},
\alpha_{2}),$ {and} $(\A, -\mu, -\Delta, \alpha_{1},
\alpha_{2})$.
\end{proposition}
\textbf{Proof:}
By a direct computation, we obtain the results. 

$ \hfill \square $
\begin{theorem}\label{relevent theorem}
Let $(\A, \mu, \Delta, (\alpha_{i})_{i=1, 2}))$ be a  finite dimensional hom-associative ternary infinitesimal bialgebra. Then so is $(\A^{\ast}, \Delta^{\ast}, \mu^{\ast}, (\alpha^{\ast}_{i})_{i=1, 2})$.
\end{theorem}
\textbf{Proof:}
One has checked that $(\A^{\ast}, \Delta^{\ast}, (\alpha^{\ast}_{i})_{i=1, 2})$ is a hom-associative tenary algebra  and
$(\A^{\ast} , \mu^{\ast}, (\alpha^{\ast}_{i})_{i=1, 2})$ is a hom-coassociative ternary coalgebra.  It remains to establish the condition 
(\ref{IATB}) for $\Delta^{\ast}$ and $\mu^{\ast}.$  For that, let us consider $a, b, c\in \A$ and $\varphi, \psi, \chi \in 
\A^{\ast}.$ We  obtain:
\beqs
\langle\mu^{\ast}\circ\Delta^{\ast}(\varphi\otimes\psi\otimes\chi), a\otimes b\otimes c\rangle &=&
\langle \varphi\otimes\psi\otimes\chi, \Delta\circ\mu(a\otimes b\otimes c)\rangle\cr
&=&\langle\varphi\otimes\psi\otimes\chi, (\mu\otimes\alpha\otimes\alpha)(\alpha\otimes\alpha\otimes\Delta)(a\otimes b\otimes c) \rangle +\cr
&&\langle \varphi\otimes\psi\otimes\chi, (\alpha\otimes\mu\otimes\alpha)(\sigma\otimes \alpha\otimes\sigma)(\alpha\otimes\Delta\otimes\alpha)(a\otimes b\otimes c)\rangle\cr
&&+\langle\varphi\otimes\psi\otimes\chi, (\alpha\otimes\alpha\otimes\mu)(\Delta\otimes\alpha\otimes\alpha)(a\otimes b\otimes c)\rangle\cr
&=&\langle(\alpha^{\ast}\otimes\alpha^{\ast}\otimes\Delta^{\ast})(\mu^{\ast}\otimes\alpha^{\ast}\otimes\alpha^{\ast})(\varphi\otimes\psi\otimes\chi), a\otimes b\otimes c\rangle +\cr
&&\langle(\alpha^{\ast}\otimes\Delta^{\ast}\otimes\alpha^{\ast})(\sigma^{\ast}\otimes\alpha^{\ast}\otimes\sigma^{\ast})
(\alpha^{\ast}\otimes\mu^{\ast}\otimes\alpha^{\ast})(\varphi\otimes\psi\otimes\chi), \cr
&& a\otimes b\otimes c\rangle + \langle(\Delta^{\ast}\otimes\alpha^{\ast}\otimes\alpha^{\ast})(\alpha^{\ast}\otimes\alpha^{\ast}\otimes\mu^{\ast})(\varphi\otimes\psi\otimes\chi),\cr
&& a\otimes b\otimes c\rangle \mbox{ where } \sigma^{\ast}(\varphi\otimes\psi)= \psi\otimes\varphi,
\eeqs
 showing that the condition (\ref{IATB}) holds for $\Delta^{\ast}$ and $\mu^{\ast}$.
$ \hfill \square $

\section{Concluding remarks}
The main results obtained in this study can be summarized as follows:
\begin{enumerate}
\item[i)] Partially and totally hom-coassociative ternary coalgebras are formulated in Definition \ref{PHCTC} and Definition \ref{THCTC}, and  characterized  in terms of the algebra constant structures in Theorem \ref{Constant structures theorem for HATC}. The relation between a hom-coassociative ternary coalgebra and its dual is given in Theorem \ref{teo-co-asso}.
\item[ii)] Theorem \ref{teo-co-asso1} connects associative $3$-ary algebras to coassociative $3$-ary coalgebras.
\item[iii)] Theorem \ref{theorem of isomorphisms and duals} links the isomorphism between two hom-coassociative ternary coalgebras to the isomorphism of their duals.
\item[iv)] The definitions of the trimodule of totally and partially hom-associative ternary algebras are given in Definition \ref{trimodule of THATA} and in Definition \ref{trimodule of PHATA},  respectively. The matched pairs of totally and partially hom-associative ternary algebras are built in Theorem \ref{Theo. Match. Pair for totally} and Theorem \ref{Theo. Match. Pair for partially}, respectively.
\item[v)] Partially and totally hom-associative ternary infinitesimal bialgebras  are defined in Definition \ref{definition of PHB} and Definition \ref{definition of HTB}, respectively,  while their dual structures are provided in Theorem \ref{pb1} and  Theorem \ref{tb1}. Finally, hom-associative ternary infinitesimal bialgebras are defined in Definition \ref{definition of inf. ternary hom-bialgebra} and the relevent proprietes are in Proposition \ref{relevent proposition}  and in Theorem \ref{relevent theorem}.
\end{enumerate}
}
 \section*{Aknowledgement}
This work is supported by TWAS Research Grant RGA No. 17 - 542 RG / MATHS / AF / AC \_G  -FR3240300147. The ICMPA-UNESCO Chair is in partnership with Daniel Iagolnitzer Foundation (DIF), France, and the Association pour la Promotion Scientifique de l'Afrique (APSA), supporting the development of mathematical physics in Africa.

 \end{document}